\documentclass[dvips]{article}

\RequirePackage[OT1]{fontenc}
\RequirePackage[ps,amsthm,amsmath,noinfoline]{imsart}
\RequirePackage{amssymb,amsfonts} 
\RequirePackage{hyperref}

\pubyear{2005}
\volume{2}
\firstpage{348}
\lastpage{384}
\doi{10.1214/154957805100000168}

\numberwithin{equation}{section}
\startlocaldefs
\theoremstyle{plain}
\newtheorem{thm}{Theorem}[section]
\newtheorem{prop}[thm]{Proposition}   
\newtheorem{cor}[thm]{Corollary}      
\newtheorem{lemma}{Lemma}[section]
\theoremstyle{definition}

\theoremstyle{remark}
\newtheorem{rem}{Remark}[section]
\endlocaldefs

\newcommand{\ds}{\displaystyle}
\newcommand{\R}{\mathbf{R}}
\newcommand{\C}{\mathbf{C}}

\newcommand{\Hy}{\mathbf{H}}
\newcommand{\F}{\mathcal{F}}
\newcommand{\Bo}{\mathcal{B}}

\newcommand{\I}{\operatorname{I}}
\newcommand{\eil}{\overset{(\text{\rm law})}{=}}
\newcommand{\pa}{\partial}
\begin{document}
\begin{frontmatter}
%
\title{Exponential functionals of \\ Brownian motion, \\
II: Some related diffusion processes\protect\thanksref{T1}}
\thankstext{T1}{This is an original survey paper.}
\runtitle{Exponential functionals of BM, II}
\runauthor{H. Matsumoto, M. Yor}

\begin{aug}
\author{\fnms{Hiroyuki} \snm{Matsumoto}\ead[label=e1]{matsu@info.human.nagoya-u.ac.jp}}
\address{Graduate School of Information Science,
Nagoya University,\\ Chikusa-ku, Nagoya 464-8601, Japan\\
\printead{e1}}

\author{\fnms{Marc} \snm{Yor}\ead[label=e2]{deaproba@proba.jussieu.fr}}
\address{Laboratoire de Probabilit\'es and Institut universitaire de France,\\
Universit\'e Pierre et Marie Curie,\\
175 rue du Chevaleret, F-75013 Paris, France\\
\printead{e2}}
\end{aug}

\begin{abstract}
This is the second part of our survey
on exponential functionals of Brownian motion.
We focus on the applications of the results
about the distributions of the exponential functionals,
which have been discussed in the first part.
Pricing formula for call options for the Asian options,
explicit expressions for the heat kernels
on hyperbolic spaces,
diffusion processes in random environments
and extensions of L\'evy's and Pitman's theorems are discussed.
\end{abstract}



\received{\smonth{9} \syear{2005}}

\end{frontmatter}
\section{Introduction}\label{s1}

Let $B=\{B_t,t\geqq0\}$ be a one-dimensional Brownian motion
starting from $0$ and defined on a probability space $(\Omega,\F,P)$.
Denoting by $B^{(\mu)}=\{B_t^{(\mu)}=B_t+\mu t\}$
the corresponding Brownian motion with constant drift $\mu\in\R$,
we consider the exponential functional
$A^{(\mu)}=\{A_t^{(\mu)}\}$ defined by
\begin{equation} \label{1e:object}
A_t^{(\mu)}=\int_0^t \exp(2B_s^{(\mu)})ds, \quad t\geqq0.
\end{equation}
In Part I \cite{my-survey1} of our survey,
we have discussed about the probability law of $A_t^{(\mu)}$
for fixed $t$ and about several related topics.

Among the results, we have shown some explicit (integral) representations
for the density of $A_t^{(\mu)}$.
In particular, we have proven the following formula
originally obtained in Yor \cite{yor-aap92}:
\begin{equation} \label{1e:yor-f}
P(A_t^{(\mu)}\in du, B_t^{(\mu)}\in dx) =
e^{\mu x-\mu^2t/2}
\exp\biggl(-\frac{1+e^{2x}}{2u}\biggr)
\theta(e^x/u,t)\frac{dudx}{u},
\end{equation}
where, for $r>0$ and $t>0$,
\begin{equation} \label{1e:theta}
\theta(r,t) = \frac{r}{(2\pi^3t)^{1/2}} e^{\pi^2/2t} \int_0^\infty
e^{-\xi^2/2t} e^{-r\cosh(\xi)} \sinh(\xi)
\sin\biggl(\frac{\pi\xi}{t}\biggr) d\xi.
\end{equation}
The function $\theta(r,t)$ appears in the representation
for the (unnormalized) density of the so-called Hartman-Watson distribution
and satisfies
\begin{equation} \label{1e:hw}
\int_0^\infty e^{-\alpha^2t/2}\theta(r,t)dt = I_{\alpha}(r),
\quad \alpha>0,
\end{equation}
where $I_\alpha$ is the usual modified Bessel function.
For details, see Part I and the references cited therein.

Another important fact, which has been used in several domains
and also discussed in Part I,
is the following identity in law due to Dufresne \cite{duf-sca90}.
Let $\mu>0$.   Then one has
\begin{equation} \label{1e:daniel}
A_\infty^{(-\mu)} \equiv \int_0^\infty \exp(2B_s^{(-\mu)})ds
\eil \frac{1}{2\gamma_\mu},
\end{equation}
where $\gamma_\mu$ is a gamma random variable with parameter $\mu$, that is,
\begin{equation*}
P(\gamma_\mu \in dx) = \frac{1}{\Gamma(\mu)} x^{\mu-1} e^{-x} dx,
\quad x\geqq0.
\end{equation*}
\indent
The purpose of this second part of our surveys is
to present some results obtained by applying the formulae and identities
mentioned in Part I to Brownian motion and some related stochastic processes.

In Section~\ref{s2} we discuss about the pricing formula
for the average option, so called, Asian option in the Black-Scholes model.

In Section~\ref{s3} we present some formulae for the heat kernels
of the semigroups generated by the Laplacians on hyperbolic spaces.
By reasoning in probabilistic terms,
we obtain not only the classical formulae
but also new expressions.

In Section~\ref{s4} we apply the results on exponential functionals
to a question pertaining to
a class of diffusion processes in random environments.

In Section~\ref{s5} we show Dufresne's recursion relation
for the probability density of $A_t^{(\mu)}$ with respect to $\mu$
which, as we have seen in Part I,
plays an important role in several studies on exponential functionals.

Dufresne's relation is important in studying extensions or analogues
of L\'evy's and Pitman's theorems
about, respectively, $\{M_t^{(\mu)}-B_t^{(\mu)}\}$ and
$\{2M_t^{(\mu)}-B_t^{(\mu)}\}$,
where $M_t^{(\mu)}=\max_{0\leqq s \leqq t}B_s^{(\mu)}$,
by means of exponential functionals.
These topics are finally discussed in Section~\ref{s6} of this survey,
where we consider the stochastic process defined by
\begin{equation*}
M_t^{(\mu),\lambda} = \frac{1}{\lambda} \log\biggl( \int_0^t
\exp(\lambda B_s^{(\mu)})ds \biggr), \quad t>0.
\end{equation*}
By the Laplace method, we easily see that, as $\lambda\to\infty$,
$M_t^{(\mu),\lambda}$ converges to $M_t^{(\mu)}$,
and we prove that $\{M_t^{(\mu),\lambda}-B_t^{(\mu)}\}$ and
$\{2M_t^{(\mu),\lambda}-B_t^{(\mu)}\}$ are diffusion processes
for any $\lambda\in\R$.
Hence, the classical L\'evy and Pitman theorems may be seen
as limiting results of those as $\lambda\to\infty$.
\section{Asian options}\label{s2}

In this section we consider the Asian or average call option
in the framework of the Black-Scholes model and
present some identities for the pricing formula.

By the Black-Scholes model, we mean a market model
which consists of a riskless bond $b=\{b_t\}$ with a constant interest rate
and a risky asset $S=\{S_t\}$
with a constant appreciation rate and volatility.
That is, letting $r>0,$ $\mu\in\R$ and $\sigma>0$ be constants,
we let $b$ and $S$ be given by the stochastic differential equation
\begin{equation*}
\frac{db_t}{b_t} = rdt, \quad
\frac{dS_t}{S_t} = \mu dt + \sigma dB_t,
\end{equation*}
where $B=\{B_t\}$ is a one-dimensional Brownian motion with $B_0=0$
defined on a complete probability space $(\Omega,\F,P)$.

For simplicity we normalize them by setting $b_0=1$.
Then we have
\begin{equation*}
b_t=\exp(rt) \quad \text{\rm and} \quad
S_t=S_0 \exp(\sigma B_t + (\mu-\sigma^2/2)t).
\end{equation*}
\indent Following the standard procedure,
we consider the discounted stock price $\widetilde{S}=\{\widetilde{S}_t\}$
given by
\begin{equation*}
\widetilde{S}_t = e^{-rt} S_t =
S_0 \exp(\sigma B_t + (\mu-r-\sigma^2/2)t).
\end{equation*}
\indent Then, by Girsanov's theorem,
there exists a unique probability measure $Q$
which is absolutely continuous with respect to $P$ and
under which $\widetilde{S}$ is a martingale.
$Q$ is called the martingale measure for $\widetilde{S}$ and we have
\begin{equation*}
\frac{dQ}{dP}\bigg|_{\F_T} = \exp\biggl( - \frac{\mu-r}{\sigma}B_T -
\frac{(\mu-r)^2}{2\sigma^2}T \biggr),
\end{equation*}
where $\F_T=\sigma\{B_s,s\leqq T\}$.
$\widetilde{B}=\{\widetilde{B}_t=B_t+\sigma^{-1}(\mu-r)t\}$ is
a Brownian motion under $Q$.

Let us consider the European and the Asian call options
with fixed strike price $k>0$ and maturity $T$.
The payoffs are given by
\begin{equation*}
(S_T-k)_{+} \quad \text{\rm and} \quad
({\mathcal A}(T)-k)_{+},
\end{equation*}
respectively, where $x_{+}=\max\{x,0\}$ and
\begin{equation*}
{\mathcal A}(t) = \frac{1}{t} \int_0^t S_u du, \quad
0 < t \leqq T.
\end{equation*}
\indent By the Black-Scholes formula or by the non-arbitrage argument,
we can show that the theoretical price $C_E(k,T)$ and $C_A(k,T)$
of these call options at time $t=0$ are given by
\begin{align*}
 & C_E(k,T) = e^{-rT} E^Q[(S_T - k)_{+}] \\
\intertext{and}
 & C_A(k,T) = e^{-rT} E^Q[({\mathcal A}(T)-k)_{+}],
\end{align*}
where $E^Q$ denotes the expectation with respect to the martingale measure $Q$.

\begin{prop} \label{op:comp}
If $r\geqq0,$ one has $C_A(k,T) \leqq C_E(k,T)$ for every $k>0$ and $T>0.$
\end{prop}

\begin{proof}
We have
\begin{equation*}
C_A(k,T) = e^{-rT} E^Q\biggl[ \biggl( \frac{1}{T} \int_0^T S_t dt - k
\biggr)_{+} \biggr].
\end{equation*}
Using Jensen's inequality, we get
\begin{equation*}
C_A(k,T) \leqq e^{-rT} \frac{1}{T} \int_0^T E^Q[(S_t-k)_{+}]dt.
\end{equation*}
\indent Since $r\geqq0$,
$\{S_t=S_0\exp(\sigma\widetilde{B}_t+(r-\sigma^2/2)t)\}$ is
a submartingale under $Q$.
Therefore, using Jensen's inequality again, we see that
$\{(S_t-k)_{+}\}$ is also a submartingale.
Hence we obtain
\begin{equation*} \begin{split}
 C_A(k,T) & \leqq e^{-rT} \frac{1}{T} \int_0^T E^Q[(S_T-k)_{+}]dt \\
 & = e^{-rT} E^Q[(S_T-k)_{+}] = C_E(k,T).
\end{split} \end{equation*}
\end{proof}

For more discussions on $C_A(k,T)$,
see Geman-Yor \cite{gy-mf93}, Rogers-Shi \cite{rs}
and the references cited therein.

By using explicit expressions for the density of $A_t^{(\mu)}$
discussed in Part I,
we obtain several integral representations for $C_A(k,T)$.
However, they are complicated.
Hence we omit this approach and
consider instead the Laplace transform of $C_A(k,T)$ in $T$.

In the following we set $\sigma=2$ and consider
\begin{equation*}
A_t^{(\mu)} = \int_0^t \exp(2B_s^{(\mu)})ds
\end{equation*}
under the original probability measure $P$
to follow the same convention as in Part I and
in other parts of the present article.

Let $T_\lambda$ be an exponential random variable
with parameter $\lambda>0$ independent of $B$.
Yor \cite{yor-cr92} (see also Part I)
has shown the identity in law
\begin{equation*}
A_{T_\lambda}^{(\mu)} \eil \frac{Z_{1,a}}{2\gamma_b},
\end{equation*}
where $a=(\nu+\mu)/2, b=(\nu-\mu)/2, \nu=\sqrt{2\lambda+\mu^2}$,
$Z_{1,a}$ is a beta variable with parameters $(1,a)$,
$\gamma_b$ is a gamma variable with parameter $b$ and
$Z_{1,a}$ and $\gamma_b$ are independent.

From this identity we deduce the following result.

\begin{thm} \label{ot:lap1}
For all $\mu\in\R,$ $\lambda>\max\{2(1+\mu),0\}$ and $k>0,$ we have
\begin{equation*} \begin{split}
\lambda \int_0^\infty e^{-\lambda t} & E[(A_t^{(\mu)}-k)_{+}] dt \\
 & = \frac{1}{(\lambda-2(1+\mu))\Gamma(b-1)} \int_0^{1/2k}
e^{-t} t^{b-2} (1-2kt)^{a+1}dt.
\end{split} \end{equation*}
\end{thm}

The same formula has been proven in \cite{gy-mf93}
with the help of some properties of Bessel processes.

We next present another proof of Theorem \ref{ot:lap1},
following Donati-Martin, Ghomrasni and Yor \cite{dgy},
who have used, as an auxiliary tool,
the stochastic process $Y^{(\mu)}(x)=\{Y_t^{(\mu)}(x)\}$ given by
\begin{equation*}
Y_t^{(\mu)}(x) = \exp(2B_t^{(\mu)}) \biggl( x + \int_0^t
\exp(-2B_s^{(\mu)}) ds \biggr).
\end{equation*}
$Y^{(\mu)}(x)$ is a diffusion process with generator
\begin{equation*}
 {\mathcal L}^{(\mu)} = 2x^2 \frac{d^2}{dx^2} + (2(1+\mu)x+1) \frac{d}{dx} =
 2x^{1-\mu} e^{1/2x} \frac{d}{dx} \biggl( x^{1+\mu} e^{-1/2x}
 \frac{d}{dx} \biggr).
\end{equation*}
In fact, in \cite{dgy}, the authors have taken advantage of
the identity in law
\begin{equation*}
\int_0^t \exp(2B_s^{(\mu)})ds \eil Y_t^{(\mu)}(0) \equiv
\exp(2B_t^{(\mu)}) \int_0^t \exp(-2B_s^{(\mu)})ds
\end{equation*}
for every fixed $t>0$
and have computed the Laplace transform of $E[(Y_t^{(\mu)}(0)-k)_{+}]$ in $t$
by using the general theory of the Sturm-Liouville operators.

We present an explicit form of
the Green function for ${\mathcal L}^{(\mu)}$.
For this purpose we recall the confluent hypergeometric functions
$\Phi(\alpha,\gamma;z)$ and $\Psi(\alpha,\gamma;z)$
of the first and second kinds defined by
\begin{align}
 & \Phi(\alpha,\gamma;z) = \sum_{k=0}^\infty \frac{(\alpha)_k}{(\gamma)_k}
 \frac{z^k}{k!} \label{oe:def-Phi} \\
\intertext{and}
 & \Psi(\alpha,\gamma;z) =
 \frac{\Gamma(1-\gamma)}{\Gamma(1+\alpha-\gamma)} \Phi(\alpha,\gamma;z) +
 \frac{\Gamma(\gamma-1)}{\Gamma(\alpha)} z^{1-\gamma}
 \Phi(1+\alpha-\gamma,2-\gamma;z), \notag
\end{align}
where $(\alpha)_0=1$ and
\begin{equation*}
(\alpha)_k=\frac{\Gamma(\alpha+k)}{\Gamma(\alpha)}=
\alpha(\alpha+1)\cdots(\alpha+k-1), \quad k=1,2,...
\end{equation*}
For details about the confluent hypergeometric functions,
we refer to Lebedev \cite{leb}.
$\Phi(\alpha,\gamma;z)$ and $\Psi(\alpha,\gamma;z)$ are
linearly independent solutions for the linear differential equation
\begin{equation*}
zu'' + (\gamma-z) u' - \alpha u = 0.
\end{equation*}
\indent We set $\nu=\sqrt{2\lambda+\mu^2}$ and define the functions
$u_1$ and $u_2$ on $(0,\infty)$ by
\begin{align}
 & u_1(x) = x^{-(\mu+\nu)/2} \Psi\biggl(\frac{\mu+\nu}{2},1+\nu; \frac{1}{2x}
 \biggr) \label{oe:u1} \\
\intertext{and}
 & u_2(x) = x^{-(\mu+\nu)/2} \Phi\biggl(\frac{\mu+\nu}{2},1+\nu; \frac{1}{2x}
 \biggr) , \label{oe:u2}
\end{align}
respectively.
Then, by straightforward computations, we can check
\begin{equation*}
{\mathcal L}^{(\mu)} u_i = \lambda u_i, \quad i=1,2.
\end{equation*}
\indent Moreover, $u_2(x)$ is monotone decreasing in $x\in(0,\infty)$.
On the other hand, recalling the integral representation
for $\Psi(\alpha,\gamma;z)$:
\begin{equation*}
\Psi(\alpha,\gamma;z) = \frac{1}{\Gamma(\alpha)} \int_0^\infty e^{-zt}
t^{\alpha-1} (1+t)^{\gamma-\alpha-1}dt, \qquad \alpha>0, z>0,
\end{equation*}
(cf. \cite{leb}, p.268),
we can easily show that $u_1(x)$ is monotone increasing.
In fact, we have
\begin{equation*}
u_1(x) = \frac{1}{\Gamma((\nu+\mu)/2)} \int_0^\infty e^{-\xi/2}
\xi^{(\mu+\nu)/2-1} (1+x\xi)^{(\nu-\mu)/2}d\xi, \quad x\geqq0.
\end{equation*}
\indent About the Wronskian, it is known (\cite{leb}, p.265) that
\begin{equation*}
\Phi(\alpha,\gamma;z) \Psi'(\alpha,\gamma;z) -
\Phi'(\alpha,\gamma;z) \Psi(\alpha,\gamma;z) =
- \frac{\Gamma(\gamma)}{\Gamma(\alpha)} z^{-\gamma} e^z,
\end{equation*}
which yields
\begin{equation*}
\frac{1}{x^{-(1+\mu)}e^{1/2x}} (u'_1(x)u_2(x)-u_1(x)u'_2(x)) =
\frac{2^\nu \Gamma(1+\nu)}{\Gamma((\mu+\nu)/2)}.
\end{equation*}
Here the function $x^{-(1+\mu)}e^{1/2x}$ is the derivative of
the scale function for $Y^{(\mu)}(x)$.

Checking the boundary conditions, we obtain the following.

\begin{prop}
Let $u_1(x)$ and $u_2(x)$ be the functions defined by
\eqref{oe:u1} and \eqref{oe:u2}{\rm.}
Then the Green function $G^{(\mu)}(x,y;\lambda)$ for ${\mathcal L}^{(\mu)}$
with respect to the Lebesgue measure is given by
\begin{equation*}
G^{(\mu)}(x,y;\lambda) = \frac{\Gamma((\mu+\nu)/2)}{2^{1+\nu}\Gamma(1+\nu)}
y^{\mu-1} e^{-1/2y} u_1(x) u_2(y), \quad 0 \leqq x \leqq y.
\end{equation*}
\end{prop}

In order to proceed to a proof of Theorem \ref{ot:lap1},
we recall the following identity presented in \cite{dgy}:
\begin{multline*}
\int_a^\infty (\xi-a) \xi^{-1-(\nu-\mu)/2} e^{-1/\xi}
\Phi\biggl(\frac{\nu+\mu}{2},1+\nu;\frac{1}{\xi}\biggr) d\xi \\
 = \frac{\Gamma((\nu-\mu)/2-1)}{\Gamma((\nu-\mu)/2+1)}
a^{1-(\nu-\mu)/2} e^{-1/a}
\Phi\biggl( \frac{\nu+\mu}{2}+2,1+\nu;\frac{1}{a} \biggr),
\end{multline*}
which may be proven by Kummer's transformation
\begin{equation*}
\Phi(\alpha,\gamma;x)=e^x\Phi(\gamma-\alpha,\gamma;-x)
\end{equation*}
and the series expansion \eqref{oe:def-Phi} of $\Phi$.
It is a special case of a general formula
given on page 279, Problem 21, Lebedev \cite{leb}.
See also \cite{gy-cr92}.

Then, noting that $u_1(0)=2^{(\mu+\nu)/2}$, we obtain
\begin{equation*} \begin{split}
\lambda & \int_0^\infty e^{-\lambda t} E[(A_t^{(\mu)}-k)_{+}] dt
 =\lambda \int_k^\infty (y-k) G^{(\mu)}(0,y;\lambda)dy \\
 & = \frac{\lambda\Gamma((\mu+\nu)/2)\Gamma((\nu-\mu)/2-1)}
 {2^{1+(\nu-\mu)/2}\Gamma(1+\nu)\Gamma((\nu-\mu)/2+1)} \\
 & \qquad \times k^{1-(\nu-\mu)/2} e^{-1/2k}
 \Phi\biggl(\frac{\mu+\nu}{2}+2,1+\nu;\frac{1}{2k}\biggr).
\end{split} \end{equation*}
\noindent Moreover, we recall the integral representation of $\Phi$:
\begin{equation*}
\Phi(\alpha,\gamma;z) =
\frac{\Gamma(\gamma)}{\Gamma(\alpha)\Gamma(\gamma-\alpha)}
\int_0^1 e^{zu} u^{\alpha-1}(1-u)^{\gamma-\alpha-1}du.
\end{equation*}
Then, after some elementary computations, we arrive at
\begin{multline*}
\lambda \int_0^\infty e^{-\lambda t} E[(A_t^{(\mu)}-k)_{+}]dt \\
 =
\frac{\lambda\Gamma((\mu+\nu)/2)}{4\Gamma((\nu-\mu)/2+1)\Gamma((\nu+\mu)/2+2)}
\int_0^{1/2k} e^{-u} b^{b-2}(1-2ku)^{a+1}du.
\end{multline*}
Finally, by using the identity $z\Gamma(z)=\Gamma(z+1)$, we obtain
\begin{equation*}
\frac{\lambda\Gamma((\mu+\nu)/2)}{4\Gamma((\nu-\mu)/2+1)\Gamma((\nu+\mu)/2+2)}
= \frac{1}{(\lambda-2(\mu+1))\Gamma((\nu-\mu)/2-1)}
\end{equation*}
and Theorem \ref{ot:lap1}.
\section{Heat kernels on hyperbolic spaces}\label{s3}

Let $\Hy^n$ be the upper half space in $\R^n$ given by
\begin{equation*}
\{z=(x,y)=(x^1,...,x^{n-1},y); x\in\R^{n-1},y>0\},
\end{equation*}
endowed with the Poincar\'e metric $ds^2=y^{-2}(dx^2+dy^2)$.
The Riemannian volume element is given by $dv=y^{-n}dxdy$ and
the distance $d(z,z')$ between $z,z'\in\Hy^n$ is given by the formula
\begin{equation} \label{he:distance}
\cosh(d(z,z')) = \frac{|x-x'|^2+y^2+(y')^2}{2yy'},
\end{equation}
where $|x-x'|$ is the Euclidean distance between $x,x'\in\R^{n-1}$.

The Laplace-Beltrami operator $\Delta_n$ is written as
\begin{equation} \label{he:laplacian}
\Delta_n = y^2 \sum_{i=1}^{n-1} \biggl(\frac{\pa}{\pa x^i}\biggr)^2 +
y^2\biggl(\frac{\pa}{\pa y}\biggr)^2 - (n-2)y\frac{\pa}{\pa y}.
\end{equation}
We denote by $p_n(t,z,z')$ the heat kernel
with respect to the volume element $dv$ of the semigroup
generated by $\Delta_n/2$.
Since $p_n(t,z,z')$ is a function of $r=d(z,z')$ for a fixed $t>0$,
we occasionally write $p_n(t,r)$ for $p_n(t,z,z')$.

Then, for $n=2$ and $3$, the following formulae are well known:
\begin{align}
 & p_2(t,r) = \frac{\sqrt{2}e^{-t/8}}{(2\pi t)^{3/2}} \int_r^\infty
 \frac{be^{-b^2/2t}}{(\cosh(b)-\cosh(r))^{1/2}}db, \label{he:cl2} \\
 & p_3(t,r) = \frac{1}{(2\pi t)^{3/2}} \frac{r}{\sinh(r)}
 \exp\biggl(-\frac{t}{2}-\frac{r^2}{2t}\biggr). \label{he:cl3}
\end{align}
Moreover, the following recursion formula due to Millson is also well known
and we also have explicit expressions of $p_n(t,r)$ for every $n\geqq 4$:
\begin{equation} \label{he:millson}
p_{n+2}(t,r) = - \frac{e^{-nt/2}}{2\pi\sinh(r)}
\frac{\pa}{\pa r}p_n(t,r).
\end{equation}
For details about the real hyperbolic space $\Hy^n$ and
the classical formulae for the heat kernels,
we refer the reader to Davies \cite{davies}.

Gruet \cite{gruet} has considered the Brownian motion on $\Hy^n$,
which is a diffusion process generated by $\Delta_n/2$, and
has derived a new integral representation for $p_n(t,r)$
by using the explicit expression \eqref{1e:yor-f}
for the joint density of $(A_t^{(\mu)},B_t^{(\mu)}).$
While the classical expressions for $p_n(t,r)$ have different forms
for odd and even dimensions,
Gruet's formula \eqref{he:gruet-f} below holds for every $n$.

\begin{thm} \label{ht:gruet}
For every $n\geqq2, t>0, z,z'\in\Hy^n,$ it holds that
\begin{equation} \label{he:gruet-f} \begin{split}
p_n(t,z,z') & = \frac{e^{-(n-1)^2t/8}}{\pi(2\pi)^{n/2}t^{1/2}}
\Gamma\biggl(\frac{n+1}{2}\biggr) \\
 & \quad \times \int_0^\infty
\frac{e^{(\pi^2-b^2)/2t}\sinh(b)\sin(\pi b/t)}
{(\cosh(b)+\cosh(r))^{(n+1)/2}}db,
\end{split} \end{equation}
where $r=d(z,z').$
\end{thm}

Before giving a proof of \eqref{he:gruet-f},
we mention its relationship to the classical formulae.
First of all we note that Millson's formula \eqref{he:millson} is
easily obtained from \eqref{he:gruet-f}
if we differentiate both hand sides of \eqref{he:gruet-f}
with respect to $r$.

When $n=3$, the integrand on the right hand side of \eqref{he:gruet-f}
may be extended to a meromorphic function in $b$ on $\C$.
Hence we can apply residue calculus and obtain \eqref{he:cl3}.

In the case $n=2$, which is the most interesting and important,
we compute the Laplace transform in $t$ of the right hand sides of
\eqref{he:cl2} and \eqref{he:gruet-f}.
Then, using the Hankel-Lipschitz formula for the modified Bessel functions
(see Watson \cite{watson}, p.386),
we can check the coincidence of the Laplace transforms or
of the expressions for the Green function.
For details, see \cite{gruet}, \cite{matsu-bull01}, \cite{mny-12}.

We give a proof of \eqref{he:gruet-f} and
see how the exponential functional $A^{(\mu)}$ comes into the story.

\smallskip

\noindent{\it Proof of Theorem \ref{ht:gruet}.}
Let $(W,\Bo,P)$ be the $n$-dimensional standard Wiener space
with the canonical filtration $\{\Bo_s\}_{s\geqq0}$:
$W$ is the space of all $\R^n$-valued continuous paths
$w_{\cdot}=(w_{\cdot}^1,...,w_{\cdot}^{n-1},w_{\cdot}^n)$ starting from $0$
with the topology of uniform convergence on compact intervals,
$\Bo$ is the topological $\sigma$-field,
$\Bo_s$ is the sub $\sigma$-field of $\Bo$
generated by $\{w_u,0\leqq u \leqq s\}$ and
$P$ is the $n$-dimensional Wiener measure.

The Brownian motion on $\Hy^n$ may be obtained as the unique solution
of the stochastic differential equation
\begin{equation*} \begin{cases}
dX^i(s)=Y(s)dw_s^i, & i=1,...,n-1, \\
dY(s)=Y(s)dw_s^n - \frac{n-2}{2}Y(s)ds. &
\end{cases} \end{equation*}
We denote by $Z_z=\{Z_z(t,w)=(X_z(t,w),Y_z(t,w)),t\geqq0\}$
the unique strong solution satisfying $Z_z(0)=z=(x,y)$.
Then we have
\begin{equation} \label{he:bm-exp} \begin{cases}
X_z^i(t,w)= x^i + \int_0^t y \exp(B_s^{(\mu)})dw_s^i, & i=1,...,n-1, \\
Y_z(t,w)=y\exp(B_t^{(\mu)}), &
\end{cases} \end{equation}
where $B_s^{(\mu)}=w_s^n+\mu s$ and $\mu=-(n-1)/2$.

Let $F_t^{(\mu)}$ be the $\R^{n-1}$-valued random variable defined by
\begin{equation*}
F_t^{(\mu)}= \biggl( \int_0^t \exp(B_s^{(\mu)})dw_s^1,...,
\int_0^t \exp(B_s^{(\mu)})dw_s^{n-1} \biggr).
\end{equation*}
Then the conditional distribution of $F_t^{(\mu)}$
given $\{Y_z(s),0\leqq s \leqq t\}$ or $\{w_s^n,0\leqq s \leqq t\}$ is
the $(n-1)$-dimensional Gaussian distribution with mean $0$ and
covariance matrix $A_t^{(\mu)}\I_{n-1}$,
$\I_{n-1}$ being the $(n-1)$-dimensional identity matrix.

Note that the heat kernel $p_n(t,z,z')$ may be written as
\begin{equation*}
p_n(t,z,z')=\int_W \widetilde{\delta}_{z'}(Z_z(t,w))dP(w) =
(y')^n \int_W \delta_{z'}(Z_z(t,w)) dP(w),
\end{equation*}
where $\widetilde{\delta}_{z'}$ and $\delta_{z'}$ are
the Dirac delta functions concentrated at $z'$
with respect to the volume element $dv$ and
the Lebesgue measure $dz=dxdy$, respectively, and
$\delta_{z'}(Z_z(t,w))$ is the composition of
the distribution $\delta_{z'}$ and the smooth Wiener functional $Z_z(t,w)$
in the sense of Malliavin calculus (see \cite{IW}).
Therefore, we obtain
\begin{equation} \label{he:to-fou} \begin{split}
 & p_n(t,z,z') \\
 & = (y')^n \int_W \delta_{(x',y')}(x+yF_t^{(\mu)},y\exp(B_t^{(\mu)}))dP \\
 & = \biggl(\frac{y'}{y}\biggr)^n \int_W \delta_{((x'-x)/y,y'/y)}
 (F_t^{(\mu)},\exp(B_t^{(\mu)}))dP \\
 & = \biggl(\frac{y'}{y}\biggr)^n \int_W
 \frac{1}{(2\pi A_t^{(\mu)})^{(n-1)/2}}
 \exp\biggl(-\frac{|x'-x|^2}{2y^2A_t^{(\mu)}}\biggr)
 \delta_{y'/y}(\exp(B_t^{(\mu)}))dP \\
 & = \biggl(\frac{y'}{y}\biggr)^{n-1} \int_W
 \frac{1}{(2\pi A_t^{(\mu)})^{(n-1)/2}}
 \exp\biggl(-\frac{|x'-x|^2}{2y^2A_t^{(\mu)}}\biggr)
 \delta_{\log(y'/y)}(B_t^{(\mu)})dP,
\end{split} \end{equation}
where we have used the same notation for the Dirac delta functions
on $\R^n$ and $\R$.

Now we apply formula \eqref{1e:yor-f} for the last expression.
Then we obtain
\begin{multline*}
p_n(t,z,z') = \biggl(\frac{y'}{y}\biggr)^{n-1} \int_0^\infty
\frac{1}{(2\pi u)^{(n-1)/2}} \exp\biggl(-\frac{|x'-x|^2}{2y^2u}\biggr) \\
 \times \biggl(\frac{y'}{y}\biggr)^{-(n-1)/2} e^{-(n-1)^2t/8}
 \frac{1}{u} \exp\biggl(-\frac{1+(y'/y)^2}{2u}\biggr)
 \theta(y'/yu,t)du.
\end{multline*}
Moreover, changing variables by $v=y'/yu$ and using \eqref{he:distance},
we obtain
\begin{equation} \label{he:to-green}
p_n(t,z,z') = \frac{e^{-(n-1)^2t/8}}{(2\pi)^{(n-1)/2}} \int_0^\infty
v^{(n-3)/2} \exp(-v\cosh(r)) \theta(v,t)dv
\end{equation}
for $r=d(z,z')$.

Finally we use the integral representation \eqref{1e:theta}
for $\theta(v,t)$.
Then, changing the order of the integrations by Fubini's theorem,
we obtain \eqref{he:gruet-f} after some elementary computations.
\hfill $\square$

\begin{rem} \label{hr:to-green}
Recalling formula \eqref{1e:hw},
we can easily obtain an explicit expression of the Green function
for $\Delta_n$ from formula \eqref{he:to-green}.
\end{rem}

\begin{rem} \label{hr:fou}
From the last expression of \eqref{he:to-fou}, we obtain
\begin{multline*}
p_n(t,z,z')=\biggl(\frac{y'}{y}\biggr)^{n-1} \int_{\R^{n-1}}
e^{-\sqrt{-1}\langle x'-x,\lambda\rangle} d\lambda \\
 \times \int_W \exp\biggl(-\frac{1}{2}|\lambda|^2y^2A_t^{(\mu)}\biggr)
\delta_{\log(y'/y)}(B_t^{(\mu)})dP.
\end{multline*}
Hence, we see that the Laplacian $\Delta_n$ on $\Hy^n$
and the Schr\"odinger operator on $\R$ with the Liouville potential
$-\frac{1}{2}\frac{d^2}{dx^2}+\frac12 |\lambda|^2e^{2x}$ are
unitary equivalent, which may be directly verified by Fourier analysis and
has been already pointed out in Comtet \cite{comtet-87},
Debiard-Gaveau \cite{dg-canada87}, Grosche \cite{grosche-88} and so on.
See also \cite{im-jfa99}.
\end{rem}

In the rest of this section, we restrict ourselves to the case $n=2$ and
consider two questions related to the results and formulae presented above.
For other related topics, see, e.g., \cite{ag-bl} and \cite{gruet-mono}.

Let us consider the following Schr\"odinger operator $H_k$, $k\in\R$,
on $\Hy^2$ with a magnetic field:
\begin{equation*}
H_k = \frac12 y^2 \biggl(\sqrt{-1}\frac{\pa}{\pa x} + \frac{k}{y}\biggr)^2
- \frac12 y^2 \biggl(\frac{\pa}{\pa y}\biggr)^2.
\end{equation*}
The differential $1$-form $\alpha=ky^{-1}dx$ is called the vector potential
and its exterior derivative $d\alpha=ky^{-2}dx\wedge dy$
represents the corresponding magnetic field.
Since $d\alpha$ is equal to constant $k$ times the volume element $dv$,
we call $H_k$ a Schr\"odinger operator with a constant magnetic field.
It is essentially the same as the Maass Laplacian
which plays an important role in several domains of mathematics,
e.g., number theory, representation theory and so on.
For details, see \cite{fay}, \cite{im-jfa99}
and the references cited therein.

In \cite{im-jfa99}, the authors have started their arguments
from the Brownian motion on $\Hy^2$ given in the above proof
of Theorem \ref{ht:gruet} and
have discussed about explicit and probabilistic expressions
for the heat kernel $q_k(t,z,z')$ of the semigroup generated by $H_k$.
They have also applied the results to a study of the Selberg trace formula
on compact quotient spaces, i.e., compact Riemannian surfaces,
and have shown close relationship between the spectrum and
the action integrals for the corresponding classical paths.
It should be mentioned that some physicists have shown similar results
in the context of Feynman path integrals prior to \cite{im-jfa99}.
See, e.g., \cite{cgo}, \cite{grosche}, \cite{GS}.

Explicit formulae for several quantities related to the operator $H_k$,
e.g., the Green functions, the heat kernels,
have been obtained by Fay \cite{fay} by harmonic analysis.
On the other hand, starting from computations by Feynman path integrals,
Comtet \cite{comtet-87} and Grosche \cite{grosche-88}
have obtained explicit forms of the Green functions.

From the point of view of probability theory along the line of
\cite{im-jfa99},
another explicit representation for the heat kernel $q_k(t,z,z')$
has been shown in \cite{ams-ln} by using an extension of
formula \eqref{1e:yor-f} and Gruet's formula \eqref{he:gruet-f}.
We introduce the result in \cite{ams-ln}
together with some arguments taken from \cite{im-jfa99}.

To show an explicit representation for $q_k(t,z,z')$,
we recall from Proposition 2.2 in \cite{im-jfa99}
(see also the references therein) that
$q_k(t,z,z')$ may be written in the form
\begin{equation} \label{he:gt}
q_k(t,z,z') = \biggl(\frac{z'-\bar{z}}{z-\bar{z'}}\biggr)^k
g_k(t,d(z,z'))
\end{equation}
for some positive function $g_k(t,r)$.
This is a consequence of the group action of $\mathrm{SL}(2;\R)$
on $\Hy^2$.
Here a point $z = (x,y) \in \Hy^2$ is identified
with $z = x+\sqrt{-1}y \in \C$,
$d(z,z')$ is the hyperbolic distance given by \eqref{he:distance}, and,
for $\omega=|\omega|\exp(\sqrt{-1}\theta)\in\C$ with $-\pi<\theta\leqq\pi$,
$\omega^k=|\omega|^k\exp(\sqrt{-1}k\theta)$.
Therefore, if $x=x'$, we have $q_k(t,z,z')=g_k(t,d(z,z'))$.

\begin{thm} \label{ht:ams}
The function $g_k(t,r)$ on the right hand side of \eqref{he:gt}
is given by
\begin{equation} \label{he:ams}
g_k(t,r) = \frac{\sqrt{2}e^{-t/8-k^2t/2}}{(2\pi t)^{3/2}} \int_r^\infty
\frac{\cosh(2k\varphi(b,r))be^{-b^2/2t}}{(\cosh(b)-\cosh(r))^{1/2}}db,
\end{equation}
where
\begin{equation*}
\varphi(b,r)=\mathrm{Argcosh}\biggl(\frac{\cosh(b/2)}{\cosh(r/2)}\biggr),
\quad 0 \leqq r \leqq b.
\end{equation*}
\end{thm}

\begin{rem}
When $k=0$, $g_0(t,r)$ coincides with the classical formula
\eqref{he:cl2} for the heat kernel $p_2(t,r)$ on $\Hy^2$.
\end{rem}

\begin{proof}
We show \eqref{he:ams} when $|k|<1/2$.
Formula \eqref{he:ams} for a general value of $k$ follows
from this result on the special case by analytic continuation.
We use the same notations as those in the proof of Theorem \ref{ht:gruet}.

Let $I_t(\alpha)$ denote the stochastic line integral (cf. \cite{IW})
of the differential $1$-form $\alpha=ky^{-1}dx$
along the path $\{Z_z(s),0\leqq s \leqq t\}$ of the Brownian motion $Z_z$
on $\Hy^2$:
\begin{equation*}
I_t(\alpha) = \int_0^t \alpha(Z_z(s)) \circ dZ_z(s).
\end{equation*}
In fact, it is easy to show
\begin{equation*}
I_t(\alpha) = \int_0^t \frac{k}{Y_z(s)}\circ dX_z(s) = kw_t^1.
\end{equation*}
\indent By using the It\^o formula, we have
\begin{equation*} \begin{split}
q_k(t,z,z') & = \int_W \exp(-\sqrt{-1}I_t(\alpha,w))
\widetilde{\delta}_{z'}(Z_z(t,w))dP(w) \\
 & = \int_W \exp(-\sqrt{-1}kw_t^1)
 \widetilde{\delta}_{(x',y')}(X_z(t,w),Y_z(t,w))dP(w).
\end{split} \end{equation*}
As in the proof of Theorem \ref{ht:gruet},
we consider the conditional distribution of
$(w_t^1,\int_0^t\exp(B_s^{(-1/2)})dw_s^1)$ given
$\{B_s^{(-1/2)}=w_s^2-s/2,0\leqq s \leqq t\}$.
Then it is easy to see that this conditional distribution is
a two-dimensional Gaussian distribution with mean $0$ and covariance matrix
\begin{equation*} \begin{pmatrix}
t & a_t^{(-1/2)} \\ a_t^{(-1/2)} & A_t^{(-1/2)}
\end{pmatrix} , \end{equation*}
where
\begin{equation*}
a_t^{(\mu)} = \int_0^t \exp(B_s^{(\mu)})ds \quad \text{\rm and} \quad
A_t^{(\mu)} = \int_0^t \exp(2B_s^{(\mu)})ds.
\end{equation*}
Taking the conditional expectation and using the Cameron-Martin theorem,
we obtain
\begin{multline*}
q_k(t,z,z') = e^{-t/8-k^2t/2} \biggl(\frac{y'}{y}\biggr)^{3/2}
\int_W \frac{1}{\sqrt{2\pi A_t(w^2)}} \\
\times \exp\biggl(-\frac{1}{2A_t(w^2)}\biggl(
\frac{x'-x}{y}+\sqrt{-1}ka_t(w^2) \biggr)^2 \biggr)
\delta_{y'/y}(\exp(w_t^2))dP,
\end{multline*}
where
\begin{equation*}
a_t(w^2) = \int_0^t \exp(w_s^2)ds \quad \text{\rm and} \quad
A_t(w^2) = \int_0^t \exp(2w_s^2)ds.
\end{equation*}
In the same way as is mentioned in Remark \ref{hr:fou}, we may write
\begin{equation} \label{he:to-morse} \begin{split}
q_k(t,z,z') & = e^{-t/8-k^2t/2} \frac{\sqrt{yy'}}{2\pi}
\int_{\R} e^{-\sqrt{-1}(x'-x)\lambda}d\lambda \\
 & \times \int_W
 \exp\biggl( -\frac12 \lambda^2y^2A_t(w^2) + \lambda k y a_t(w^2) \biggr)
 \delta_{\log(y'/y)}(w_t^2)dP \\
 & = e^{-t/8-k^2t/2} \frac{\sqrt{yy'}}{2\pi} \int_{\R}
 e^{-\sqrt{-1}(x'-x)\lambda} q_{\lambda,k}(t,\log y,\log y') d\lambda,
\end{split} \end{equation}
where $q_{\lambda,k}(t,\xi,\eta)$ denotes the heat kernel of the semigroup
generated by the Schr\"odinger operator $H_{\lambda,k}$ on $\R$
with the Morse potential given by
\begin{equation*}
H_{\lambda,k}=-\frac12 \frac{d^2}{d\xi^2} + V_{\lambda,k}, \qquad
V_{\lambda,k}(\xi)=\frac12 \lambda^2e^{2\xi} - \lambda k e^{\xi}.
\end{equation*}
\indent In \cite{ams-ln} and Part I,
we have shown an explicit representation
for $q_{\lambda,k}(t,\xi,\eta)$: for $\lambda>0$,
\begin{equation} \label{he:morse}
q_{\lambda,k}(t,\xi,\eta) = \int_0^\infty e^{2ku} \frac{1}{2\sinh(u)}
\exp(-\lambda(e^{\xi}+e^{\eta})\coth(u))
\theta(\phi,t/4)du,
\end{equation}
where the function $\theta(r,t)$ is given by \eqref{1e:theta} and
$\phi=2\lambda e^{(\xi+\eta)/2}/\sinh(u)$.

For $\lambda<0$, we have
$q_{\lambda,k}(t,\xi,\eta)=q_{-\lambda,-k}(t,\xi,\eta)$.

We now recall the remark following the statement of Theorem \ref{ht:ams}
and consider the case $x'=x$.
Then, combining \eqref{he:to-morse} and \eqref{he:morse}, we obtain
\begin{equation*} \begin{split}
 & q_k(t,z,z') = g_k(t,r) \\
 & = e^{-t/8-k^2t/2} \frac{\sqrt{yy'}}{2\pi} \int_0^\infty
 (q_{\lambda,k}(t,\log y,\log y')+q_{-\lambda,k}(t,\log y,\log y'))d\lambda \\
 & = e^{-t/8-k-2t/2} \frac{\sqrt{yy'}}{2\pi} \int_0^\infty d\lambda
 \int_0^\infty \frac{\cosh(2ku)}{\sinh(u)}
 \exp(-\lambda(y+y')\coth(u)) \\
 & \hspace{7cm} \times
 \theta\biggl(\frac{2\lambda\sqrt{yy'}}{\sinh(u)},\frac{t}{4} \biggr)du.
\end{split} \end{equation*}
Note that the integral is convergent if $|k|<1/2$.

Then, using the integral representation \eqref{1e:theta}
for $\theta(r,t)$ and carrying out the integration in $\lambda$ first,
we obtain
\begin{align*}
 & g_k(t,r) = \frac{e^{-t/8-k^2t/2}}{\pi(2\pi)^{3/2}t^{1/2}}
 \int_0^\infty \cosh(2ku)F_t(u)du, \\
 & F_t(u) = \int_0^\infty
 \frac{e^{2(\pi^2-\xi^2)/t}\sinh(\xi)\sin(4\pi\xi/t)}
 {(\cosh(r/2)\cosh(u)+\cosh(\xi))^2}d\xi.
\end{align*}
By Gruet's formula \eqref{he:gruet-f}, we have
\begin{equation*} \begin{split}
p_3(t/4,r) & = \frac{2e^{-t/8}}{\pi(2\pi)^{3/2}t} \int_0^\infty
\frac{e^{2(\pi^2-b^2)/t}\sinh(b)\sin(4\pi b/t)}
{(\cosh(b)+\cosh(r))^2}db \\
 & = \biggl(\frac{2}{\pi t}\biggr)^{3/2} \frac{r}{\sinh(r)}
 \exp\biggl(-\frac{t}{8}-\frac{r^2}{t}\biggr).
\end{split} \end{equation*}
It is now easy to show \eqref{he:ams} from these formulae.
\end{proof}

Similar arguments to those in the proofs of
Theorems \ref{ht:gruet} and \ref{ht:ams} are available
to study the Laplace-Beltrami operators on the complex and quaternion
hyperbolic spaces.
Also on these symmetric spaces of rank one,
we have explicit expressions of Brownian motions as Wiener functionals and
we can show explicit representations for the heat kernels
and for the Green functions.
For details, see \cite{matsu-bull01}.

Next we consider the diffusion process on $\Hy^2$
associated to the infinitesimal generator
\begin{equation*}
{\mathcal L}_{\nu,\mu} = \frac12 y^2 \biggl(\frac{\pa}{\pa x}\biggr)^2 +
\frac12 y^2 \biggl(\frac{\pa}{\pa y}\biggr)^2 - \nu y \frac{\pa}{\pa x} -
\biggl(\mu-\frac12 \biggr) y \frac{\pa}{\pa y},
\end{equation*}
where $\nu\geqq0$ and $\mu>0$.
The operator ${\mathcal L}_{\nu,\mu}$ is invariant
under the special transforms
on $\Hy^2$ of the form $z\mapsto az+b,$ $a>0$ and $b\in\R$,
while the operator $H_k$ and, in particular, the Laplacian $\Delta_2$ are
invariant under the action of $\mathrm{SL}(2;\R)$.

The diffusion process starting from $z=(x,y)$
associated to $\mathcal{L}_{\nu,\mu}$
may be realized as the unique strong solution
$\{Z_t^{(\nu,\mu)}=(X_t^{(\nu,\mu)},Y_t^{(\nu,\mu)})\}$
of the stochastic differential equation
\begin{equation*} \begin{cases}
\ds dX_t = Y_t dw_t^1 - \nu Y_t dt, & X_0=x, \\
\ds dY_t = Y_t dw_t^2 - \biggl(\mu-\frac12 \biggr) Y_t dt, & Y_0=y,
\end{cases} \end{equation*}
defined on a two-dimensional Wiener space.
As in the case of Brownian motion on $\Hy^2$,
$Z^{(\nu,\mu)}$ is also represented as a Wiener functional by
\begin{equation*} \begin{cases}
\ds X_t^{(\nu,\mu)}= x + y \int_0^t \exp(B_s^{(-\mu)})dW_s^{(-\nu)}, \\
\ds Y_t^{(\nu,\mu)}= y \exp(B_t^{(-\mu)}),
\end{cases} \end{equation*}
where $W_s^{(-\nu)}=w_s^1-\nu s$ and $B_s^{(-\mu)}=w_s^2-\mu s.$

Since $\mu$ is assumed to be positive,
$Y_t^{(\nu,\mu)}$ converges to $0$ as $t$ tends to $\infty$.
Following \cite{bcfy}, we show that the distribution of $X_t^{(\nu,\mu)}$
converges weakly as $t\to\infty$
and that we can specify the limiting distribution.
It is enough to consider the special case $x=0$ and $y=1$.

\begin{thm} \label{ht:bcfy}
When $x=0$ and $y=1,$
the distribution of $X_t^{(\nu,\mu)}$ on $\R$ converges weakly
to the distribution with density
\begin{equation*}
f(\xi) = C_{\nu,\mu}
\frac{\exp(-2\nu \mathrm{Arctan}(\xi))}{(1+\xi^2)^{\mu+1/2}},
\quad \xi\in\R.
\end{equation*}
\end{thm}

For details on the normalizing constant $C_{\nu,\mu}$, see \cite{bcfy}.
Note that, if $\nu=0$ and $\mu=1/2$, that is,
if we consider the Brownian motion on $\Hy^2$,
the limiting distribution is the Cauchy distribution
as in the case of the hitting distribution on lines
of standard Brownian motion on $\R^2$.
In general, the limiting distribution belongs to
the type IV family of Pearson distributions (cf. \cite{jkb}).

It should be mentioned that the functional
$\int_0^\infty \exp(B_s^{(-\mu)})dW_s^{(-\nu)}$ has been much studied
in the context of risk theory.
See Paulsen \cite{paulsen} and the references cited therein about this.
In \cite{paulsen} the density is derived when $\nu>1$.
See also \cite{ag-bl} and \cite{bcf} about some results
in special cases.
For further related discussions, see \cite{my-osaka} and \cite{yor-meander}.

We present a probabilistic proof taken from \cite{bcfy},
where we also find an analytic proof.

\begin{proof}
The limiting distribution coincides with that for the stochastic process
$\{\bar{X}_t^{(\nu,\mu)}\}$ given by
\begin{equation*}
\bar{X}_t^{(\nu,\mu)} = x \exp(B_t^{(-\mu)}) +
\int_0^t \exp(B_s^{(-\mu)})dW_s^{(-\nu)}.
\end{equation*}
We also consider the diffusion process $\{\widetilde{X}_t^{(\nu,\mu)}\}$
given by
\begin{equation*}
\widetilde{X}_t^{(\nu,\mu)} = \exp(B_t^{(-\mu)}) \biggl( x +
\int_0^t \exp(-B_s^{(-\mu)})dW_s^{(-\nu)} \biggr)
\end{equation*}
with infinitesimal generator
\begin{equation*}
\widetilde{{\mathcal L}}^{(\nu,\mu)} = \frac{1+x^2}{2} \frac{d^2}{dx^2} -
\biggl( \nu + \biggl( \mu - \frac12 \biggr)x \biggr) \frac{d}{dx}.
\end{equation*}

By the invariance of the law of Brownian motion
under time reversal from a fixed time,
$\bar{X}_t^{(\nu,\mu)}$ and $\widetilde{X}_t^{(\nu,\mu)}$ are
identical in law for any fixed $t>0$.
Therefore, to prove the theorem,
we only have to check $\widetilde{{\mathcal L}}^{(\nu,\mu)*}f=0$
for the adjoint operator $\widetilde{{\mathcal L}}^{(\nu,\mu)*}$
to $\widetilde{\mathcal{L}}^{(\nu,\mu)}$.
\end{proof}

\begin{rem}
Set
\begin{equation*}
A_t^{(-\mu)} = \int_0^t \exp(2B_s^{(-\mu)})ds \quad \text{\rm and} \quad
a_t^{(-\mu)} = \int_0^t \exp(B_s^{(-\mu)})ds.
\end{equation*}
The joint distribution of $(A_t^{(-\mu)},a_t^{(-\mu)},B_t^{(-\mu)})$ or
the Laplace transform of the conditional distribution of $A_t^{(-\mu)}$
given $(a_t^{(-\mu)},B_t^{(-\mu)})$ has been studied
in \cite{ams-ln} (see also Part I).

We have, using some obvious notation,
\begin{equation*}
X_\infty^{(\nu,\mu)} \eil \gamma_{A_\infty^{(-\mu)}} -
\nu a_\infty^{(-\mu)}
\end{equation*}
for a Brownian motion $\{\gamma_t\}$ independent of $B$.
Hence we obtain
\begin{equation*}
f(\xi) = E\biggl[ \frac{1}{\sqrt{2\pi A_\infty^{(-\mu)}}}
\exp\biggl(-\frac{(\xi+\nu a_\infty^{(-\mu)})^2}{2A_\infty^{(-\mu)}}
\biggr) \biggr].
\end{equation*}
However, we have not succeeded in obtaining Theorem \ref{ht:bcfy}
from this expression.
\end{rem}

\begin{rem}
The limiting distribution with density $f(\xi)$ belongs to
the domain of attraction of a stable distribution,
whose characteristic function $\phi$ is of the form
\begin{align*}
 & \phi(t) = \exp\biggl
(\sqrt{-1}zt+c|t|^\alpha\biggl(1+\sqrt{-1}\gamma\mathrm{sgn}(t)
\tan\biggl(\frac{\alpha\pi}{2}\biggr)\biggr)\biggr), \quad
0<\alpha\leqq 2, \alpha\ne1, \\
\intertext{or}
 & \phi(t)=\exp\biggl(\sqrt{-1}zt+c|t|\biggl(1\sqrt{-1}\gamma\mathrm{sgn}(t)
\log(|t|) \frac{2}{\pi} \biggr) \biggr),
\end{align*}
where $c>0, -1<\gamma<1$ and $z\in\R$.

It is also the case if we consider the hitting distribution
on $\{\mathrm{Im}(z)=a\}$, that is,
the distribution of $X_{\tau_a}^{(\nu,\mu)}$ when
$Y_0^{(\nu,\mu)} = y > a > 0,$
where $\tau_a$ is the first hitting time of $a$ by $\{Y_t^{(\nu,\mu)}\}$.
For details, see the original paper \cite{bcfy}.
\end{rem}
\section{Maximum of a diffusion process in random environment}\label{s4}

The purpose of this section is to survey the work by Kawazu-Tanaka \cite{kt}
on the maximum of a diffusion process in a drifted random environment.
In \cite{kt},
several equalities and inequalities for the exponential functionals
of Brownian motion are used.

Let $W=\{W(y),y\in\R\}$ be a Brownian environment
defined on a probability space $(\Omega_1,\F_1,P_1)$:
$\{W(y),y\geqq0\}$ and $\{W(-y),y\geqq0\}$ are
independent one-dimensional Brownian motions with $W(0)=0$.
For $c\in\R$, we set $W^{(c)}(y)=W(y)+cy$.

For each $\omega_1\in\Omega_1$,
we consider a diffusion process
$X(W(\omega_1))=X(W)=\{X(t,W),t\geqq0\}$, $W(\omega_1)=W(\cdot, \omega_1)$,
with $X(0,W)=0$, whose infinitesimal generator is given by
\begin{equation*}
\frac12 \exp(W^{(c)}(x)) \frac{d}{dx} \biggl( \exp(-W^{(c)}(x))
\frac{d}{dx} \biggr) .
\end{equation*}
We denote by $P_{\omega_1}$ the probability law of the diffusion process
$\{X(t,W(\omega_1))\}$ and consider the full probability law
${\mathcal P} = \int P_{\omega_1} P_1(d\omega_1)$ of $\{X(t,\cdot)\}$.

A scale function $S^{(c)}(x)=S_W^{(c)}(x)$ for $X(W)$ is given by
\begin{equation*}
S^{(c)}(x)=\int_0^x \exp(W^{(c)}(y))dy \ \text{\rm for} \ x \geqq0, \quad
= - \int_x^0 \exp(W^{(c)}(y))dy \ \text{\rm for} \ x\leqq0.
\end{equation*}
By the general theory of the one-dimensional diffusion processes,
$\{S^{(c)}(X_t(W))\}$ may be represented
as a random time change of another Brownian motion
and, based on this representation,
several interesting results have been obtained.
For these results, see, e.g., Brox \cite{brox},
Hu-Shi-Yor \cite{hsy}, Kawazu-Tanaka \cite{kt-2}.

In the rest of this section, we assume $c>0$.
Then we have $S^{(c)}(\infty)=\infty$ and $S^{(c)}(-\infty)>-\infty$,
and therefore, $\max_{t\geqq0}X(t)<\infty$, ${\mathcal P}$-a.s.

The question we discuss in the present section is how the tail probability
${\mathcal P}(\max_{t\geqq0}X(t)>x)$ decays as $x\to\infty$.
We have
\begin{equation} \label{rme:basis}
{\mathcal P}\left(\max_{t\geqq0}X(t)>x\right) =
E^{P_1}\biggl[\frac{-S^{(c)}(-\infty)}{S^{(c)}(x)-S^{(c)}(-\infty)}\biggr],
\quad x>0,
\end{equation}
and several results on the exponential functional given by \eqref{1e:object}
are quite useful in this study.
The random variables $S^{(c)}(x)$ and $S^{(c)}(-\infty)$ are independent.
We also note (cf. \eqref{1e:daniel})
that $-S^{(c)}(-\infty)$ is distributed as $2\gamma_{2c}^{-1}$,
where $\gamma_{2c}$ is a gamma random variable with parameter $2c$.

\begin{thm} \label{rmt:kt}
{\rm (i)} If $c>1,$ then one has
\begin{equation*}
{\mathcal P}\left(\max_{t\geqq0}X(t)>x\right) = \frac{2(c-1)}{2c-1}
\exp\biggl(-\biggl(c-\frac{1}{2}\biggr)x\biggr) (1+o(1)),
\quad x\to\infty.
\end{equation*}
{\rm (ii)} If $c=1,$ then
\begin{equation*}
{\mathcal P}\left(\max_{t\geqq0}X(t) > x\right) =
\sqrt{\frac{2}{\pi}} x^{-1/2} e^{-x/2} (1+o(1)), \quad x\to\infty.
\end{equation*}
{\rm (iii)} If $0<c<1,$ then
\begin{equation*}
{\mathcal P}\left(\max_{t\geqq0}X(t)>x\right) =
C x^{-3/2} e^{-c^2x/2} (1+o(1)),
\quad x\to\infty,
\end{equation*}
where the constant $C$ is given by
\begin{align*}
 & C= \frac{2^{5/2-2c}}{\Gamma(2c)} \int_0^\infty \cdots \int_0^\infty
 \frac{za^{2c-1}e^{-a/2}}{a+z}y^{2c} e^{-\lambda(y,u)z} u \sinh(u) \;
 dadydzdu, \\
 & \lambda(y,u) = \frac{1+y^2}{2} + y \cosh(u).
\end{align*}
\end{thm}

Before proceeding to a proof for each assertion,
we rewrite the right hand side of \eqref{rme:basis} into different forms.
We set $A^{(c)}=-S^{(c)}(-\infty)$ and
\begin{equation*}
f^{(c)}(a,x) = E^{P_1}[(a + S^{(c)}(x))^{-1}], \qquad a>0, x>0.
\end{equation*}
Then we have
\begin{align}
{\mathcal P}\left(\max_{t\geqq0}X(t)>x\right) & = E[A^{(c)}
f^{(c)}(A^{(c)},x)] \label{rme:1} \\
 & = \frac{1}{\Gamma(2c)} \int_0^\infty \frac{2}{\xi}
 f^{(c)}\biggl(\frac{2}{\xi},x\biggr) \xi^{2c-1} e^{-\xi} d\xi. \label{rme:2}
\end{align}
Moreover, considering the time reversal of $W$, we easily obtain
\begin{multline*}
f^{(c)}(a,x) = e^{-(c-1/2)x} E^{P_1}\biggl[ \biggl(
a \exp(W^{(-c)}(x)) + \int_0^x \exp(W^{(-c)}(y))dy \biggr)^{-1} \\
 \times \exp(W(x)-x/2) \biggr].
\end{multline*}
By the Cameron-Martin theorem, we also obtain
\begin{equation} \label{rme:for1}
f^{(c)}(a,x) = e^{-(c-1/2)x} E^{P_1}\biggl[ \biggl( a \exp(W^{(1-c)}(x)) +
\int_0^x \exp(W^{(1-c)}(y))dy \biggr)^{-1} \biggr].
\end{equation}
By considering the time reversal again, we may write
\begin{equation*}
f^{(c)}(a,x) = e^{-(c-1/2)x} E^{P_1}\biggl[ \biggl( a +
\int_0^x \exp(W^{(c-1)}(y)) dy \biggr)^{-1}
\exp(W^{(c-1)}(x)) \biggr].
\end{equation*}

\noindent{\it Proof of} (i).
From \eqref{rme:1} and \eqref{rme:for1}, we have
\begin{equation*} \begin{split}
e^{(c-1/2)x} & {\mathcal P}\left(\max_{t\geqq0}X(t)>x\right) \\
 & = E^{P_1}\biggl[ A^{(c)} \biggl( A^{(c)} \exp(W^{(1-c)}(x)) +
\int_0^x \exp(W^{(1-c)}(y)) dy \biggr)^{-1} \biggr]
\end{split} \end{equation*}
and, by the independence of $A^{(c)}$ and $\{W^{(1-c)}(y),y\geqq0\}$,
we obtain
\begin{equation*}
\lim_{x\to\infty} e^{(c-1/2)x} {\mathcal P}\left(\max_{t\geqq0}X(t)>x\right)
= E\biggl[\frac{2}{\gamma_{2c}}\biggr]
E\biggl[ \frac{\gamma_{2(c-1)}}{2} \biggr],
\end{equation*}
where $\gamma_{\mu}$ is a gamma random variable with parameter $\mu>0$.
Easy evaluation of the right hand side yields the assertion.
\hfill $\square$

Before proceeding to a proof of (ii), we prepare two lemmas.

\begin{lemma} \label{rml:1}
Setting
\begin{equation*}
\psi(x) = E^{P_1}\biggl[ \biggl(\int_0^x \exp(W(y)) dy \biggr)^{-1}
\exp(W(x)) \biggr],
\end{equation*}
we have
\begin{equation*}
\psi(x) = E^{P_1}\biggl[ \biggl( \int_0^x \exp(W(y))dy \biggr)^{-1} \biggr]
\end{equation*}
and
\begin{equation} \label{rme:psi}
\lim_{x\to\infty} \sqrt{2\pi x}\psi(x) = 1.
\end{equation}
\end{lemma}

\begin{proof}
The first assertion is easily shown by time reversal.

We can show the second assertion from the identity
\begin{equation*}
E^{P_1}\biggl[ \biggl( \int_0^x \exp(2W(y))dy \biggr)^{-1} \bigg|
W(x)=u \biggl] =
\frac{ue^{-u}}{x\sinh(u)}, \quad x>0, u\in\R,
\end{equation*}
which has been shown in Part I, Proposition 5.9.
However, we give another direct proof.

Set
\begin{equation*}
\varphi(x) = E^{P_1}\biggl[ \log\biggl(\int_0^x \exp(W(y))dy\biggr)\biggr].
\end{equation*}
Then we have $\varphi'(x)=\psi(x)$ and, if we show
\begin{equation*}
\lim_{x\to\infty} \frac{\varphi(x)}{\sqrt{x}} = \sqrt{\frac{2}{\pi}},
\end{equation*}
we obtain \eqref{rme:psi} by L'Hospital's theorem.

By the scaling property of Brownian motion, we have
\begin{equation*}
\frac{\varphi(x)}{\sqrt{x}} = E^{P_1}\biggl[ \frac{1}{\sqrt{x}}
\log\biggl( \int_0^1 \exp(\sqrt{x}W(y)) dy \biggr) \biggr] +
\frac{1}{\sqrt{x}} \log(x).
\end{equation*}
and, by the Laplace principle, we also have
\begin{equation*}
\lim_{x\to\infty} \frac{1}{\sqrt{x}} \log \biggl( \int_0^1
\exp(\sqrt{x}W(y)) dy \biggr) = \max_{0\leqq y \leqq 1}W(y).
\end{equation*}
Hence, applying the dominated convergence theorem, we obtain
\begin{equation*}
\lim_{x\to\infty}\frac{\varphi(x)}{\sqrt{x}} =
E^{P_1}\left[ \max_{0\leqq y \leqq 1}W(y)\right] = \sqrt{\frac{2}{\pi}}.
\end{equation*}
\end{proof}

\begin{lemma} \label{rml:2}
For all $x>0,$ one has
\begin{equation*}
E^{P_1}\biggl[ \biggl( \int_0^x \exp(W(y)) dy \biggr)^{-2}
\exp(W(x)) \biggr] \leqq \biggl\{ \psi\biggl(\frac{x}{2}\biggr)\biggr\}^{1/2}.
\end{equation*}
\end{lemma}

\begin{proof} We write
\begin{align*}
 & E^{P_1}\biggl[ \biggl( \int_0^x \exp(W(y)) dy \biggr)^{-2}
\exp(W(x)) \biggr] \\
 & \leqq E^{P_1}\biggl[ \biggl( \int_0^{x/2} \exp(W(y))dy \biggr)^{-1}
 \biggl( \int_{x/2}^x \exp(W(y))dy \biggr)^{-1} \exp(W(x)) \biggr] \\
 & = E^{P_1}\biggl[ \biggl( \int_0^{x/2} \exp(W(y)) dy \biggr)^{-1} \\
 & \qquad \quad \times
 \biggl( \int_{x/2}^x \exp(W(y)-W(x/2)) dy \biggr)^{-1}
 \exp(W(x)-W(x/2)) \biggr].
\end{align*}
Hence, using the independence of increments of Brownian motion, we obtain
\begin{align*}
 & E^{P_1}\biggl[ \biggl( \int_0^x \exp(W(y)) dy \biggr)^{-2}
\exp(W(x)) \biggr] \\
 & \leqq E^{P_1}\biggl[ \biggl( \int_0^{x/2} \exp(W(y)) dy \biggr)^{-1}\biggr]
 E^{P_1}\biggl[ \biggl( \int_{0}^{x/2} \exp(W(y)) dy \biggr)^{-1}
 \exp(W(x/2)) \biggr] \\
 & = \{\psi(x/2)\}^2.
\end{align*}
\end{proof}

\noindent{\it Proof of} (ii).
Set $A^{(1)}=-S^{(1)}(-\infty)$.
Then, since $A^{(1)}$ and $B_x^{(1)}$ are independent, we have
\begin{equation*} \begin{split}
E^{P_1}\biggl[ \frac{A^{(1)}}{S^{(1)}(x)}\biggr] & =
E^{P_1}[A^{(1)}] E^{P_1}\biggl[ \frac{1}{S^{(1)}(x)}\biggr] \\
 & = E\biggl[\frac{2}{\gamma_2}\biggr] E^{P_1}\biggl[
 \biggl( \int_0^x \exp(W(y))dy \biggr)^{-1} \exp(W(x)-x/2) \biggr],
\end{split} \end{equation*}
where $\gamma_2$ is a gamma variable with parameter $2$ and
we have used the Cameron-Martin theorem for the second equality.
Therefore we obtain
\begin{equation} \label{rme:si}
\lim_{x\to\infty} x^{1/2}e^{x/2}
E^{P_1}\biggl[ \frac{A^{(1)}}{S^{(1)}(x)}\biggr]
= 2 \lim_{x\to\infty}\sqrt{x}\psi(x) = \sqrt{\frac{2}{\pi}}
\end{equation}
from \eqref{rme:psi}.

We next prove
\begin{equation} \label{rme:ti}
E^{P_1}\biggl[ \frac{A^{(1)}}{S^{(1)}(x)} \biggr] -
E^{P_1}\biggl[ \frac{A^{(1)}}{A^{(1)}+S^{(1)}(x)} \biggr] \leqq
C x^{-3/4} e^{-x/2}
\end{equation}
for some absolute constant $C$.
Combining this with \eqref{rme:si} above, we obtain the assertion.

For this purpose we note the elementary inequality
\begin{equation*}
0 \leqq \frac{a}{b} - \frac{a}{a+b} \leqq
\frac{1}{2} \biggl(\frac{a}{b}\biggr)^{3/2}, \quad a,b>0.
\end{equation*}
Then we obtain
\begin{equation*}
E^{P_1}\biggl[ \frac{A^{(1)}}{S^{(1)}(x)} \biggr] -
E^{P_1}\biggl[ \frac{A^{(1)}}{A^{(1)}+S^{(1)}(x)} \biggr] \leqq
\frac12 E^{P_1}[(A^{(1)})^{3/2}] E[(S^{(1)}(x))^{-3/2}].
\end{equation*}
\indent For the first term on the right hand side, we have
\begin{equation*}
E^{P_1}[(A^{(1)})^{3/2}] = 2^{3/2} \int_0^\infty x^{-3/2} x e^{-x}dx
< \infty.
\end{equation*}
For the second term, we use the Cauchy-Schwarz inequality to show
\begin{align*}
 & E^{P_1}[(S^{(1)}(x))^{-3/2}] = E^{P_1}\biggl[ \biggl( \int_0^x
 \exp(W(y))dy \biggr)^{-3/2} \exp(W(x)-x/2) \biggr] \\
 & \leqq e^{-x/2} \bigg\{ E^{P_1}\biggl[ \biggl( \int_0^x
 \exp(W(y)) dy \biggr)^{-1} \exp(W(x)) \biggr] \biggr\}^{1/2} \\
 & \qquad \quad \qquad \times
 \biggl\{ E^{P_1}\biggl[ \biggl( \int_0^x \exp(W(y))dy \biggr)^{-2}
 \exp(W(x)) \biggr] \biggr\}^{1/2}.
\end{align*}
Then, using Lemmas \ref{rml:1} and \ref{rml:2},
we obtain \eqref{rme:ti} and the result of (ii).  \hfill $\square$

\medskip

\noindent{\it Proof of} (iii).
We prove this case by using formula \eqref{1e:yor-f}.
To do this in a direct way, we note that
\begin{equation*}
S^{(c)}_x = \int_0^x \exp(W^{(c)}(y))dy \eil
4 \int_0^{x/4} \exp(2W^{(2c)}(y)) dy,
\end{equation*}
and that the latter is $4A_{x/4}^{(2c)}$,
where $A_t^{(\mu)}$ is defined by \eqref{1e:object}.

Then, by using \eqref{rme:basis} and \eqref{1e:yor-f}, we obtain
\begin{align*}
 & {\mathcal P}\left(\max_{t\geqq0}X(t)>x\right) =
 E^{P_1}[A^{(c)} f^{(c)}(A^{(c)},x)] \\
 & = \int_0^\infty \frac{2}{a} \frac{1}{\Gamma(2c)} a^{2c-1} e^{-a} da
 \cdot e^{-c^2x/2} \int_0^\infty du \int_{\R} d\xi \int_0^\infty db \\
 & \quad \times
 e^{2c\xi} u^{-1} \exp\biggl(-\frac{1+e^{2\xi}}{2u}\biggr)
 \biggl(\frac{2}{a} + 4u \biggr)^{-1} \\
 & \quad \times
 \sqrt{\frac{2}{\pi^3x}} \frac{e^\xi}{u} e^{(2\pi^2-2b^2)/x}
 \exp\biggl(-\frac{e^\xi\cosh(b)}{u}\biggr) \sinh(b)
 \sin\biggl(\frac{4\pi\xi}{x}\biggr).
\end{align*}
From this identity we see that
the order of decay is $x^{-3/2}e^{-c^2x/2}$ and,
by using the dominated convergence theorem and
changing the variables in the integration,
we obtain the assertion.
For details, see the original paper \cite{kt}. \hfill $\square$
\section{Exponential functionals with different drifts}\label{s5}

The purpose of this section is to show a relationship
between the laws of the exponential functionals of Brownian motions
with different drifts.

In this and the next sections,
we consider several stochastic processes or transforms on path space
related to the exponential functional $\{A_t^{(\mu)}\}$.
In particular, the following transform $Z$ plays an important role.
For a continuous function $\phi:[0,\infty)\to\R$,
we associate $A(\phi)=\{A_t(\phi)\}$ and $Z(\phi)=\{Z_t(\phi)\}$ defined by
\begin{equation} \label{de:a-and-z}
 A_t(\phi) = \int_0^t \exp(2\phi(s))ds
\quad \text{ and } \quad
 Z_t(\phi) = \exp(-\phi(t)) A_t(\phi).
\end{equation}
\indent Let $\nu<\mu$ and consider two exponential functionals
$A^{(\nu)}=A(B^{(\nu)})$ and $A^{(\mu)}=A(B^{(\mu)})$:
\begin{equation*}
A_t^{(\nu)} = \int_0^t \exp(2B_s^{(\nu)})ds
\quad \text{\rm and} \quad
A_t^{(\mu)} = \int_0^t \exp(2B_s^{(\mu)})ds,
\end{equation*}
where $B_s^{(\nu)}=B_s+\nu s$ and $B=\{B_s\}$ is a one-dimensional
Brownian motion with $B_0=0$ as in the previous sections.

We first consider the case where $\nu=-\mu$ and $\mu>0$ and,
using the result in this special case,
we will give the general result in Theorem \ref{dt:general} below.

\begin{thm} \label{dt:osaka}
Let $\mu>0$ and let $\gamma_\mu$ be a gamma random variable
with density $\Gamma(\mu)^{-1}x^{\mu-1}e^{-x},$ $x>0,$
independent of $B.$
Then one has the identity in law
\begin{equation} \label{de:osaka}
\biggl\{ \frac{1}{A_t^{(-\mu)}}, t>0 \biggr\} \eil
\biggl\{ \frac{1}{A_t^{(\mu)}} + 2\gamma_\mu, t>0 \biggr\}.
\end{equation}
\end{thm}

The identity in law for a fixed $t>0$ has been obtained
by Dufresne \cite{duf-osaka} and
the extension \eqref{de:osaka} at the process level has been
given in \cite{my-osaka}.

\begin{proof}
We sketch a proof based on the theory of initial enlargements
of filtrations (cf. Yor \cite{yor-eth2}).
Another proof based on some properties of Bessel processes has been given
in \cite{my-osaka}.

Let $\mathcal{B}_t^{(-\mu)}=\sigma\{B_s^{(-\mu)},s\leqq t\}$ and
$\widehat{\mathcal{B}}_t^{(-\mu)}=
\mathcal{B}_t^{(-\mu)}\vee \sigma\{A_\infty^{(-\mu)}\}.$
Then we can show that there exists a
$\{\widehat{\mathcal{B}}_t^{(-\mu)}\}$-Brownian motion
$\{\widehat{B}_t\}$ independent of $A_\infty^{(-\mu)}$ such that
\begin{equation} \label{de:enlarge}
B_t^{(-\mu)} = \widehat{B}_t + \mu t -
\int_0^t \frac{\exp(2B_s^{(-\mu)})}{A_\infty^{(-\mu)}-A_s^{(-\mu)}}ds.
\end{equation}
We regard this identity as an equation for $\{B_s^{(-\mu)}\}$
with given initial data $\{\widehat{B}_t\}$ and $A_\infty^{(-\mu)}$.
Then, solving \eqref{de:enlarge}, we obtain
\begin{equation*}
B_t^{(-\mu)} = \widehat{B}_t^{(\mu)} -
\log\biggl( 1 + \frac{\widehat{A}_t^{(\mu)}}{A_\infty^{(-\mu)}} \biggr),
\end{equation*}
where $\widehat{B}_t^{(\mu)}=\widehat{B}_t + \mu t$ and
$\{\widehat{A}_t^{(\mu)}\}=\{A_t(\widehat{B}^{(\mu)})\}$.

As a consequence, it follows that
\begin{equation*}
A_t^{(-\mu)} = \int_0^t
\biggl( 1 + \frac{\widehat{A}_t^{(\mu)}}{A_\infty^{(-\mu)}} \biggr)^{-2}
\exp(2\widehat{B}_s^{(\mu)})ds =
\frac{\widehat{A}_t^{(\mu)}}{1+\widehat{A}_t^{(\mu)}/A_\infty^{(-\mu)}},
\end{equation*}
hence
\begin{equation*}
\frac{1}{A_t^{(-\mu)}} = \frac{1}{\widehat{A}_t^{(\mu)}} +
\frac{1}{A_\infty^{(-\mu)}}.
\end{equation*}
Since $(A_\infty^{(-\mu)})^{-1}\eil 2\gamma_\mu$ by \eqref{1e:daniel} and
$\{\widehat{B}_t^{(\mu)}\}$ is independent of $A_\infty^{(-\mu)}$,
we obtain the theorem.
\end{proof}

Next we consider the stochastic processes
$Z^{(-\mu)}=Z(B^{(-\mu)})$ and $\widehat{Z}^{(\mu)}=Z(\widehat{B}^{(\mu)})$:
\begin{equation*}
Z_t^{(-\mu)} = \exp(-B_t^{(-\mu)}) A_t^{(-\mu)}
\quad \text{\rm and} \quad
\widehat{Z}_t^{(\mu)} = \exp(-\widehat{B}_t^{(\mu)})
\widehat{A}_t^{(\mu)}.
\end{equation*}
Then we have
\begin{equation} \label{de:diff}
\frac{d}{dt}\biggl(\frac{1}{A_t^{(-\mu)}}\biggr) =
\frac{1}{(Z_t^{(-\mu)})^2}.
\end{equation}
Hence, from the identity \eqref{de:osaka},
we obtain $Z(B^{(-\mu)})\eil Z(B^{(\mu)})$ for any $\mu>0$.
Moreover, by \eqref{de:enlarge}, we have
the pathwise identity $Z_t^{(-\mu)}=\widehat{Z}_t^{(\mu)}$, $t\geqq0$.

The study of the stochastic process $Z(B^{(\mu)})$ is in fact
the main object of the next section.
We will show that, for any $\mu\in\R$,
$Z(B^{(\mu)})=\{Z_t^{(\mu)}\}$ is a diffusion process
with respect to its natural filtration $\{\mathcal{Z}_t^{(\mu)}\}$,
$\mathcal{Z}_t^{(\mu)}=\sigma\{Z_s^{(\mu)},s\leqq t\}$,
and that this result gives rise to an extension of Pitman's theorem
(\cite{pitman},\cite{rp}).

A key fact in the proof of the above mentioned result is
the following Proposition \ref{dp:key},
which also plays an important role in the rest of this section.

Before mentioning the proposition, we note another important fact.
By \eqref{de:diff}, we easily obtain the following: for every $t>0$,
\begin{equation} \label{de:filt}
\mathcal{B}_t^{(\mu)} = \sigma\{A_s^{(\mu)},s\leqq t\} =
\mathcal{Z}_t^{(\mu)} \vee \sigma\{A_t^{(\mu)}\}.
\end{equation}
In particular, $\mathcal{Z}_t^{(\mu)}$ is strictly smaller
than the original filtration $\mathcal{B}_t^{(\mu)}$
of the Brownian motion $B$,
as is also shown very clearly in the next proposition.

\begin{prop} \label{dp:key}
Let $\mu\in\R.$
Then{\rm ,} for any $t>0,$ the conditional distribution of
$\exp(B_t^{(\mu)})$ given $\mathcal{Z}_t^{(\mu)}$ is
a generalized inverse Gaussian distribution and is given by
\begin{equation} \label{de:gig}
P(\exp(B_t^{(\mu)}) \in dx | \mathcal{Z}_t^{(\mu)} ) =
\frac{1}{2K_\mu(1/z)} x^{\mu-1}
\exp\biggl(-\frac{1}{2z} \biggl( x+\frac{1}{x} \biggr)\biggr) dx, \ x>0,
\end{equation}
where $z=Z_t^{(\mu)}$ and
$K_\mu$ is the modified Bessel (Macdonald) function{\rm .}
\end{prop}

\begin{cor}
Let $\mu=0$ and write $B_t,$ $\mathcal{Z}_t$
for $B_t^{(0)},$ $\mathcal{Z}^{(0)},$ respectively{\rm.}
Then{\rm ,}\\
{\rm (i)} the conditional distributions of $\exp(B_t)$ and $\exp(-B_t)$
given $\mathcal{Z}_t$ are identical{\rm;} \\
{\rm (ii)} letting $\gamma_\delta$ be a gamma random variable
independent of $\{B_t\},$ one has
\begin{equation} \label{de:cor}
E[ f(e^{B_t}) e^{\delta B_t} | \mathcal{Z}_t ] =
E[ f(e^{B_t}+2z\gamma_\delta) e^{-\delta B_t} | \mathcal{Z}_t].
\end{equation}
\end{cor}

\noindent{\it Proof of the corollary}.\
The first assertion is easily obtained.

For the second assertion,
we consider random variables $I_z^{(\pm\delta)}$
whose densities are given by the right hand side of \eqref{de:gig},
replacing $\pm\delta$ for $\mu$.
Then, assuming that $I_z^{(-\delta)}$ and $\gamma_z$ are independent,
we have $I^{(\delta)}_z\eil I^{(-\delta)}_z+2z\gamma_\delta$.
In fact, more general identities in law
for generalized inverse Gaussian and gamma random variables are well known.
See Seshadri \cite{seshadri}, \cite{my-nmj01} and the references therein.

Hence, from the identity \eqref{de:gig} considered for $\mu=0$, we obtain
\begin{equation*} \begin{split}
E[ f(e^{B_t}) e^{\delta B_t} | \mathcal{Z}_t ] & =
\frac{K_\delta(1/z)}{K_0(1/z)} E[f(I^{(\delta)}_z)] \\
 & = \frac{K_\delta(1/z)}{K_0(1/z)}
 E[f(I^{(-\delta)}_z+2z\gamma_\delta)] \\
 & = E[ f(e^{B_t}+2z\gamma_\delta) e^{-\delta B_t} | \mathcal{Z}_t].
\end{split} \end{equation*}
\hfill $\square$

We postpone a proof of Proposition \ref{dp:key} to the next section and,
admitting this proposition as proven,
we show a general relationship between the probability laws of
the exponential functionals of Brownian motions with different drifts.

\begin{thm} \label{dt:general}
Let $\nu<\mu$ and set $\delta=(\mu-\nu)/2$ and $m=(\mu+\nu)/2.$
Then{\rm ,} for every $t\geqq0$ and for every non-negative functional $F$
on $C([0,t]\to\R),$ one has
\begin{equation} \label{de:general} \begin{split}
e^{\nu^2t/2} & E\biggl[ F\biggl( \frac{1}{A_s^{(\nu)}}, s \leqq t \biggr)
\biggl(\frac{1}{A_t^{(\nu)}}\biggr)^m \biggl] \\
 & = e^{\mu^2t/2} E\biggl[ F\biggl( \frac{1}{A_s^{(\mu)}} + 2\gamma_\delta,
s \leqq t \biggr) \biggl(\frac{1}{A_t^{(\mu)}}\biggr)^m \biggl],
\end{split} \end{equation}
where $\gamma_\delta$ is a gamma random variable with parameter $\delta$
independent of $\{B_s^{(\mu)}\}.$
\end{thm}

\begin{proof}
We start from \eqref{de:cor}.
Then, for any non-negative function $\psi$ on $\R_{+}$, we have
\begin{equation*}
E\biggl[ \psi\biggl(\frac{e^{B_t}}{z}\biggr) e^{\delta B_t}
 \bigg| \mathcal{Z}_t \biggr] =
E\biggl[ \psi\biggl(\frac{e^{B_t}}{z} + 2 \gamma_\delta\biggr)
e^{-\delta B_t} \bigg| \mathcal{Z}_t \biggr].
\end{equation*}
We also deduce from the first assertion of the corollary
\begin{equation*}
E\biggl[ \psi\biggl(\frac{e^{-B_t}}{z}\biggr)e^{-\delta B_t}
\bigg| \mathcal{Z}_t \biggr] =
E\biggl[ \psi\biggl(\frac{e^{-B_t}}{z} + 2 \gamma_\delta\biggr)
e^{\delta B_t} \bigg| \mathcal{Z}_t \biggr]
\end{equation*}
or, since $Z_t=\exp(-B_t)A_t$,
\begin{equation*}
E\biggl[ \psi\biggl(\frac{1}{A_t}\biggr) \biggl(\frac{z}{A_t}\biggr)^\delta
\bigg| \mathcal{Z}_t \biggr] =
E\biggl[ \psi\biggl(\frac{1}{A_t} + 2 \gamma_\delta\biggr)
\biggl(\frac{z}{A_t}\biggr)^{-\delta} \bigg| \mathcal{Z}_t \biggr].
\end{equation*}
Multiplying by $(Z_t)^{-\nu}$ on both hand sides,
we rewrite the last identity into
\begin{equation*}
E\biggl[ \psi\biggl(\frac{1}{A_t}\biggr) \biggl(\frac{1}{A_t}\biggr)^m
e^{\nu B_t} \bigg| \mathcal{Z}_t \biggr] =
E\biggl[ \psi\biggl(\frac{1}{A_t} + 2 \gamma_\delta\biggr)
\biggl(\frac{1}{A_t}\biggr)^m e^{\mu B_t} \bigg| \mathcal{Z}_t \biggr].
\end{equation*}
Then we obtain, for any non-negative functional $\widetilde{F}$,
\begin{equation*}
E\biggl[ \widetilde{F}(Z_s, s\leqq t) \psi\biggl(\frac{1}{A_t}\biggr)
\frac{\exp(\nu B_t)}{(A_t)^m} \biggr] =
E\biggl[ \widetilde{F}(Z_s, s\leqq t)
\psi\biggl(\frac{1}{A_t} + 2 \gamma_\delta\biggr)
\frac{\exp(\mu B_t)}{(A_t)^m} \biggr] .
\end{equation*}
\indent By the Cameron-Martin theorem, we now obtain
\begin{equation*} \begin{split}
e^{\nu^2t/2} & E\biggl[ \widetilde{F}(Z_s^{(\nu)}, s\leqq t)
\psi\biggl(\frac{1}{A_t^{(\nu)}}\biggr)
\biggl(\frac{1}{A_t^{(\nu)}}\biggr)^m \biggr] \\
 & = e^{\mu^2t/2} E\biggl[ \widetilde{F}(Z_s^{(\mu)}, s\leqq t)
\psi\biggl(\frac{1}{A_t^{(\mu)}} + 2 \gamma_\delta\biggr)
\biggl(\frac{1}{A_t^{(\mu)}}\biggr)^m \biggr].
\end{split} \end{equation*}
This identity is equivalent to \eqref{de:general}
because of \eqref{de:diff} or \eqref{de:filt}.
\end{proof}

Let $f^{(\mu)}(a,t)$ be the density of $(2A_t^{(\mu)})^{-1}$.
Then the following ``recursion'' formula,
originally due to Dufresne \cite{duf-aap01},
is deduced from \eqref{de:general}.

\begin{thm} \label{dt:density}
Let $\nu<\mu.$
Then{\rm ,} for any $t>0,$ one has{\rm ,} with $\delta=(\mu-\nu)/2,$
\begin{equation*}
e^{\nu^2t/2} f^{(\nu)}(a,t) = e^{\mu^2t/2}
\frac{a^{-m}e^{-a}}{\Gamma(\delta)} \int_0^a (a-b)^{\delta-1} b^m e^b
f^{(\mu)}(b,t)db, \quad a>0.
\end{equation*}
\end{thm}
\section{Some exponential analogues of L\'evy's and Pitman's theorems}\label{s6}

In this section we consider the two stochastic processes
$\xi^{(\mu)}=\{\xi_t^{(\mu)}\}$ and $Z^{(\mu)}=\{Z^{(\mu)}_t\}$ defined by
\begin{align*}
 & \xi^{(\mu)}_t = \exp(-2B^{(\mu)}_t) A^{(\mu)}_t =
 \exp(-2B^{(\mu)}_t) \int_0^t \exp(2B^{(\mu)}_s)ds \\
\intertext{and}
 & Z^{(\mu)}_t = \exp(-B^{(\mu)}_t) A^{(\mu)}_t,
\end{align*}
where $B^{(\mu)}_t=B_t+\mu t$ and
$B=\{B_t\}$ is a one-dimensional Brownian motion starting from $0$.

Our purpose here is to show that both $\xi^{(\mu)}$ and $Z^{(\mu)}$ are
diffusion processes, that is,
they give representations for some diffusion processes starting from $0$,
with respect to their natural filtrations and
that this result gives rise to analogues or extensions of
the celebrated L\'evy and Pitman theorems.

We start by recalling these classical theorems.
Set
\begin{equation*}
M^{(\mu)}_t = \max_{0\leqq s \leqq t}B^{(\mu)}_s.
\end{equation*}
Then the L\'evy and Pitman theorems may be stated
in the following general form with any $\mu\in\R$.

\begin{thm} \label{lpt:levy}
{\rm (i)} Let $X^{(\mu)}=\{X^{(\mu)}_t\}$ be
the bang-bang process with $X^{(\mu)}_0=0$ and with parameter $\mu,$
that is{\rm ,} the diffusion process with infinitesimal generator
$\frac12 \frac{d^2}{dx^2}-\mu \mathrm{sgn}(x)\frac{d}{dx}$ and
let $\{\ell^{(\mu)}_t\}$ be the local time of $X^{(\mu)}$ at $0.$
Then{\rm ,} the following identity in law holds {\rm :}
\begin{equation*}
\{(M^{(\mu)}_t-B^{(\mu)}_t,M^{(\mu)}_t),t\geqq0\} \eil
\{(|X^{(\mu)}_t|,\ell^{(\mu)}_t),t\geqq0\}.
\end{equation*}
{\rm (ii)} $\sigma\{M^{(\mu)}_s-B^{(\mu)}_s,s\leqq t\} =
\mathcal{B}^{(\mu)}_t \equiv \sigma\{B^{(\mu)}_s,s\leqq t\}.$
\end{thm}

\begin{thm} \label{lpt:pitman}
{\rm (i)} Let $\{\rho^{(\mu)}_t\}$ be the diffusion process
starting from $0$ with infinitesimal generator
\begin{equation*}
\frac12 \frac{d^2}{dx^2} + \mu \coth(\mu x)\frac{d}{dx}
\end{equation*}
and set $j^{(\mu)}_t=\inf_{s\geqq t}\rho^{(\mu)}_s.$
Then{\rm ,} the following identity in law holds {\rm :}
\begin{equation*}
\{(2M^{(\mu)}_t-B^{(\mu)}_t,M^{(\mu)}_t), t\geqq0\} \eil
\{(\rho^{(\mu)}_t,j^{(\mu)}_t), t\geqq0\}.
\end{equation*}
{\rm (ii)} $\sigma\{2M^{(\mu)}_s-B^{(\mu)}_s,s\leqq t\}
\subsetneqq \mathcal{B}^{(\mu)}_t.$ \\
{\rm (iii)} As a consequence of {\rm (i),}
the diffusion processes $\{2M^{(-\mu)}_t-B^{(-\mu)}_t\}$ and
$\{2M^{(\mu)}_t-B^{(\mu)}_t\}$ have the same probability law{\rm .}
\end{thm}

When $\mu=0$, $\{|X^{(0)}_t|\}$ and $\{\rho^{(0)}_t\}$ are respectively
a reflecting Brownian motion and a three-dimensional Bessel process,
and the theorems give their representations
in terms of the maximum process $\{M^{(0)}_t\}$
of a Brownian motion $B$.
For the proofs and related references,
see \cite{my-cr99}, \cite{my-nmj01}, \cite{ry} {\it et al}.
It should be noted that
the stochastic processes $\{kM^{(\mu)}_t-B^{(\mu)}_t\}$, $k\in\R$, are
not Markovian
except for these two interesting cases $k=1,2$ and the trivial case $k=0$.
For a rigorous and detailed proof, see \cite{mo}.

Now, for $\lambda>0$, we set
\begin{equation*}
M^{(\mu),\lambda}_t = \frac{1}{2\lambda} \log\biggl(
\int_0^t \exp(2\lambda B^{(\mu)}_s) ds \biggr).
\end{equation*}
Then, the Laplace principle implies
\begin{equation*}
\lim_{\lambda\to\infty} M^{(\mu),\lambda}_t = M^{(\mu)}_t,
\end{equation*}
and, by the scaling property of Brownian motion, we deduce
\begin{equation*}
\{M^{(\mu),\lambda}_t,t>0\} \eil
\biggl\{ \frac{1}{\lambda} M^{(\mu/\lambda),1}_{\lambda^2 t} -
\frac{1}{2\lambda} \log \lambda^2, t>0 \biggr\}.
\end{equation*}
Moreover, we have
\begin{align*}
 & \log(\xi^{(\mu)}_t) = \log\biggl( \int_0^t \exp(2B^{(\mu)}_s)ds \biggr)
 - 2B^{(\mu)}_t \equiv 2(M^{(\mu),1}_t - B^{(\mu)}_t) \\
\intertext{and}
 & \log(Z^{(\mu)}_t) = \log\biggl( \int_0^t \exp(2B^{(\mu)}_s)ds \biggr)
 - B^{(\mu)}_t \equiv 2M^{(\mu),1}_t - B^{(\mu)}_t.
\end{align*}
Hence, if we show that $\xi^{(\mu)}$ and $Z^{(\mu)}$ are diffusion processes
for every $\mu\in\R$, then we see that
$\{M^{(\mu),\lambda}_t-B^{(\mu)}_t\}$ and
$\{2M^{(\mu),\lambda}_t-B^{(\mu)}_t\}$ are also diffusion processes
for every $\lambda>0$ and
we can recover the L\'evy and Pitman theorems as the limiting cases
by letting $\lambda\to\infty$.

Furthermore, there is an exponential analogue of the second assertion of
Theorem \ref{lpt:pitman}
since we have shown in the previous section (see \eqref{de:filt}) that
$\mathcal{Z}^{(\mu)}_t \equiv \sigma\{Z^{(\mu)}_s,s\leqq t\}
\subsetneqq \mathcal{B}^{(\mu)}_t$ holds for every $t>0$.

We can easily show that $\xi^{(\mu)}$ is a diffusion process.
In fact, from the It\^o formula, we deduce
\begin{equation*}
d\xi^{(\mu)}_t = -2\xi^{(\mu)}_t dB_t +
((2-2\mu)\xi^{(\mu)}_t + 1) dt,
\end{equation*}
which implies the following.

\begin{thm} \label{lpt:xi}
Let $\mu\in\R.$   \\
{\rm (i)} $\xi^{(\mu)}$ is a diffusion process with respect to
the natural filtration $\{\mathcal{B}^{(\mu)}_t\}$ of $B^{(\mu)}$
and its infinitesimal generator is given by
\begin{equation*}
2x^2 \frac{d^2}{dx^2} + ((2-2\mu)x+1) \frac{d}{dx}.
\end{equation*}
{\rm (ii)} For every fixed $t>0,$ one has
\begin{equation*}
\xi^{(\mu)}_t \eil \int_0^t \exp(-2B^{(\mu)}_s)ds \eil
\int_0^t \exp(2B_s^{(-\mu)})ds.
\end{equation*}
\end{thm}

On the other hand,
it is not as easy to show that $Z^{(\mu)}$ is a diffusion process.
By the It\^o formula, we have
\begin{equation} \label{lpe:ito-z}
Z^{(\mu)}_t = - \int_0^t Z^{(\mu)}_s dB_s +
\int_0^t \biggl( \frac12 - \mu \biggr) Z^{(\mu)}_s ds +
\int_0^t \exp(B^{(\mu)}_s)ds
\end{equation}
and we need to take care of the third term on the right hand side.

Here we recall Proposition \ref{dp:key}, which implies
\begin{equation*}
E[\exp(B^{(\mu)}_s) | \mathcal{Z}^{(\mu)}_s ] =
\biggl(\frac{K_{1+\mu}}{K_\mu}\biggr)\biggl(\frac{1}{z}\biggr),
\qquad z=Z^{(\mu)}_s,
\end{equation*}
by the integral representation for $K_\mu$ (cf. \cite{leb}, p.119)
\begin{equation} \label{lpe:int-k}
K_\mu(y) = \frac12 \biggl(\frac{y}{2}\biggr)^\mu \int_0^\infty
e^{-t-y^2/4t} t^{-\mu-1}dt.
\end{equation}

Then, by replacing in \eqref{lpe:ito-z} $\exp(B_s^{(\mu)})$
by its projection on $\mathcal{Z}^{(\mu)}_s$
(cf., e.g., Liptser-Shiryaev \cite{ls}, Theorem 7.17),
we see that there exists a $\{\mathcal{Z}^{(\mu)}_t\}$-Brownian motion
$\{\beta_t\}$ such that
\begin{equation*}
Z^{(\mu)}_t = \int_0^t Z^{(\mu)}_s d\beta_s +
\biggl( \frac12 - \mu \biggr) \int_0^t Z^{(\mu)}_s ds +
\int_0^t \biggl( \frac{K_{1+\mu}}{K_\mu} \biggr)
\biggl(\frac{1}{Z^{(\mu)}_s}\biggr) ds.
\end{equation*}
Hence, admitting Proposition \ref{dp:key} as proven,
we have obtained the following.

\begin{thm} \label{lpt:zett}
Let $\mu\in\R.$  \\
{\rm (i)} $Z^{(\mu)}$ is a diffusion process on $[0,\infty)$
with respect to its natural filtration $\{\mathcal{Z}^{(\mu)}_t\}$
whose infinitesimal generator is given by
\begin{equation*}
\frac12 z^2 \frac{d^2}{dz^2} + \biggl\{ \biggl( \frac12 -\mu \biggr)z +
\biggl( \frac{K_{1+\mu}}{K_\mu} \biggr) \biggl(\frac{1}{z}\biggr)
\biggr\} \frac{d}{dz}.
\end{equation*}
{\rm (ii)} For every $t>0,$
$\mathcal{Z}^{(\mu)}_t \subsetneqq \mathcal{B}^{(\mu)}_t.$ \\
{\rm (iii)} The diffusion processes $Z^{(-\mu)}$ and $Z^{(\mu)}$ have
the same probability law{\rm .}
\end{thm}

To compare with the original L\'evy and Pitman theorems,
we present the following, which can be obtained from
Theorems \ref{lpt:xi} and \ref{lpt:zett}
by the scaling property of Brownian motion.
We call a diffusion process a Brownian motion with drift $b$
if its generator is given by $\frac12 \frac{d^2}{dx^2}+b(x)\frac{d}{dx}$.

\begin{thm}
{\rm (i)} For any $\lambda>0,$ the stochastic process
\begin{equation*}
\biggl\{\frac{1}{2\lambda} \log\biggl( \int_0^t \exp(2\lambda B_s^{(\mu)})
ds \biggr) - B^{(\mu)}_t + \frac{1}{2\lambda} \log\lambda^2, t>0
\biggr\}
\end{equation*}
is a Brownian motion with drift
\begin{equation*}
c^{(\mu),\lambda}(x) = -\mu + \frac12 \lambda e^{-2\lambda x}.
\end{equation*}
{\rm (ii)} The stochastic process
\begin{equation*}
\biggl\{\frac{1}{\lambda} \log\biggl( \int_0^t \exp(2\lambda B_s^{(\mu)})
ds \biggr) - B^{(\mu)}_t + \lambda \log\lambda^2, t>0 \biggr\}
\end{equation*}
is a Brownian motion with drift
\begin{equation*}
b^{(\mu),\lambda}(x) = -\mu + \lambda e^{-\lambda x}
\biggl(\frac{K_{1+\mu/\lambda}}{K_{\mu/\lambda}}\biggr)(e^{-\lambda x}) =
\frac{d}{dx}\bigl( \log K_{\mu/\lambda}(e^{-\lambda x})\bigr).
\end{equation*}
\end{thm}

\begin{rem}
By using the integral representation for $K_\mu(x)$
\begin{equation*}
K_\mu(x) = \frac{2^{\mu}\Gamma(\mu+1/2)}{x^\mu\sqrt{\pi}}
\int_0^\infty \frac{\cos(xt)}{(1+t^2)^{\mu+1/2}}dt, \quad x>0, \mu>0,
\end{equation*}
(cf. \cite{leb}, p.140), we can show
\begin{equation*}
\lim_{\lambda\to\infty} b^{(\mu),\lambda}(x) = \mu \coth(\mu x).
\end{equation*}
\end{rem}

The rest of this section is devoted to a proof of Proposition \ref{dp:key},
which has played an important role
not only in this section but also in the previous section
and in Part I of our survey.
To show the proposition, we prove an identity
for an anticipative transform on path space,
which may be regarded as an example of the Ramer-Kusuoka formula.
For the Ramer-Kusuoka formula,
see \cite{bf}, \cite{kusuoka}, \cite{ramer} and \cite{yano}.

Another proof for the proposition
which uses several properties of Bessel processes has been given
in \cite{my-nmj00} and
a proof based on Theorem \ref{dt:osaka},
featuring the generalized Gaussian inverse distributions, has been given
in \cite{my-nmj01}.

For our purpose we consider one more transform on path space.
For an $\R$-valued continuous function $\phi$ on $[0,\infty)$,
we define $T_\alpha(\phi)=\{T_\alpha(\phi)_t,t\geqq0\},$ $\alpha\geqq0,$ by
\begin{equation} \label{lpe:t-trans}
T_\alpha(\phi)_t = \phi(t) - \log(1+\alpha A_t(\phi)).
\end{equation}
\indent We now summarize some properties of these transforms on path space
which are easy to prove but play important roles in the following.

\begin{prop} \label{lpp:easy}
Letting $A$ and $Z$ be the transforms on path space
defined by \eqref{de:a-and-z} and
$T$ be defined by \eqref{lpe:t-trans}{\rm,} one obtains
\begin{equation*} \begin{split}
 & {\rm (i)} \frac{1}{A_t(T_\alpha(\phi))} = \frac{1}{A_t(\phi)}+\alpha,
 \qquad \qquad {\rm (ii)} Z \circ T_\alpha=Z, \\
 & {\rm (iii)} T_\alpha \circ T_\beta =T_{\alpha+\beta}, \quad \alpha,\beta>0.
\end{split} \end{equation*}
\end{prop}

The next theorem gives an example of the Ramer-Kusuoka formula.
On the left hand side of \eqref{lpe:rk} below,
the transform $T_{\alpha/e_t^{(\mu)}}$ depends on $e_t^{(\mu)}$ and
it is natural to call it anticipative.
For more discussions about this transform,
related topics and references, see \cite{dmy-rims02}.

\begin{thm} \label{lpt:rk}
Let $\mu\in\R,\alpha\geqq0$ and let $F$ be
a non-negative functional on the path space $C([0,t]\to\R).$
Then{\rm,} setting $e_t^{(\mu)}=\exp(B_t^{(\mu)}),$ one has
\begin{equation} \label{lpe:rk}
E[F(T_{\alpha/e_t^{(\mu)}}(B^{(\mu)})_s, s\leqq t)] =
E[F(B^{(\mu)}_s,s\leqq t) \Gamma_\alpha^{(\mu)}(e_t^{(\mu)},Z_t^{(\mu)})]
\end{equation}
for every $t>0,$ where
$\Gamma_\alpha^{(\mu)}(x,z)=(1+\alpha z)^\mu \Gamma_\alpha(x,z)$ and
\begin{equation*}
\Gamma_\alpha(x,z) = \exp\biggl( -\frac{\alpha}{2} \biggl( x -
\frac{1}{(1+\alpha z)x} \biggr) \biggr).
\end{equation*}
\end{thm}

\begin{proof}
We divide our proof into three steps.   \\
Step 1.\  First, we prove that
\begin{equation} \label{lpe:41}
E\biggl[\exp(-\eta e_t) F\biggl(\frac{1}{A_s}, s\leqq t\biggr) \biggr] =
E\biggl[\exp(-\eta/e_t) F\biggl(\frac{1}{A_s}+2\eta/e_t, s\leqq t
\biggr) \biggr]
\end{equation}
holds for every $\eta>0$, where $e_t=e_t^{(0)}$.
From Proposition \ref{lpp:easy},
we see that this identity is equivalent to
\begin{equation} \label{lpe:42}
E[\exp(-\eta e_t) F(B_s, s\leqq t)] = E[\exp(-\eta/e_t)
F(T_{2\eta/e_t}(B)_s, s\leqq t)].
\end{equation}
Step 2.\  We prove \eqref{lpe:rk} with $\mu=0$
from \eqref{lpe:42}.   \\
Step 3.\  We prove \eqref{lpe:rk} for general values of $\mu$.

We start from Theorem \ref{dt:osaka}, which says that
\begin{equation} \label{lpe:osaka}
E\biggl[ F\biggl( \frac{1}{A_s^{(\mu)}}+2\gamma_\mu, s\leqq t \biggr)\biggr]
= E\biggl[ F\biggl(\frac{1}{A_s^{(-\mu)}}, s\leqq t\biggr) \biggr]
\end{equation}
holds for any $\mu>0$ and
for any non-negative functional $F$ on $C([0,t]\to\R)$,
where $\gamma_\mu$ is a gamma random variable
with parameter $\mu$ independent of $B^{(\mu)}$.

We rewrite the left hand side of \eqref{lpe:osaka} into the following way:
\begin{align}
 & E\biggl[ F\biggl(\frac{1}{A_s^{(\mu)}}+2\gamma_\mu,s\leqq t\biggr)
 \biggr] \notag \\
 & = \frac{1}{\Gamma(\mu)} \int_0^\infty \eta^{\mu-1} e^{-\eta}
 E\biggl[ F\biggl( \frac{1}{A_s^{(\mu)}}+2\eta, s\leqq t\biggr)\biggr]
 d\eta \notag \\
 & = \frac{1}{\Gamma(\mu)} \int_0^\infty \eta^{\mu-1} e^{-\eta}
 E\biggl[ F\biggl( \frac{1}{A_s}+2\eta, s\leqq t\biggr)
 \exp(\mu B_t-\mu^2t/2) \biggr]  d\eta, \notag \\
 & = \frac{e^{-\mu^2t/2}}{\Gamma(\mu)} \int_0^\infty e^{-\eta}
 E\biggl[ F\biggl( \frac{1}{A_s} + 2\eta , s\leqq t\biggr) (\eta e_t)^\mu
 \biggr]  \frac{d\eta}{\eta} , \label{lpe:a}
\end{align}
where we have used the Cameron-Martin theorem for the second identity.
For the right hand side of \eqref{lpe:osaka}, we rewrite
\begin{align}
 & E\biggl[ F\biggl( \frac{1}{A_s^{(-\mu)}}, s\leqq t \biggr) \biggr] \notag \\
 & = E\biggl[ F\biggl( \frac{1}{A_s}, s\leqq t \biggr)
 \exp(-\mu B_t-\mu^2t/2) \biggr] \notag \\
 & = \frac{e^{-\mu^2t/2}}{\Gamma(\mu)} \int_0^\infty \eta^{\mu-1} e^{-\eta}
 E\biggl[F\biggl(\frac{1}{A_s}, s\leqq t\biggr) \exp(-\mu B_t) \biggr]
 d\eta \notag \\
 & = \frac{e^{-\mu^2t/2}}{\Gamma(\mu)} \int_0^\infty e^{-\eta}
 E\biggl[ F\biggl( \frac{1}{A_s}, s\leqq t\biggr)
 (\eta/e_t)^\mu \biggr]  \frac{d\eta}{\eta}.  \label{lpe:b}
\end{align}
Now, comparing \eqref{lpe:a} and \eqref{lpe:b}, we get
\begin{equation*} \begin{split}
\int_0^\infty e^{-\eta} & E\biggl[
F\biggl( \frac{1}{A_s}+2\eta, s\leqq t\biggr)
(\eta e_t)^\mu \biggr] \frac{d\eta}{\eta} \\
 & = \int_0^\infty e^{-\eta} E\biggl[ F\biggl( \frac{1}{A_s}, s\leqq t\biggr)
(\eta/e_t)^\mu \biggr] \frac{d\eta}{\eta}.
\end{split} \end{equation*}
\indent Since this identity holds for any $\mu>0$, we obtain
\begin{equation*} \begin{split}
\int_0^\infty e^{-\eta} & E\biggl[
F\biggl( \frac{1}{A_s}+2\eta, s\leqq t\biggr)
f(\eta e_t) \biggr] \frac{d\eta}{\eta} \\
 & = \int_0^\infty e^{-\eta} E\biggl[ F\biggl( \frac{1}{A_s}, s\leqq t\biggr)
f(\eta/e_t) \biggr] \frac{d\eta}{\eta}
\end{split} \end{equation*}
for any non-negative Borel function $f$ on $[0,\infty)$.
Since the left hand side is equal to
\begin{equation*}
\int_0^\infty E[\exp(-\eta/e_t) F\biggl(\frac{1}{A_s}+2\eta/e_t,
s\leqq t\biggr) \biggr] f(\eta) \frac{d\eta}{\eta}
\end{equation*}
and the right hand side is
\begin{equation*}
\int_0^\infty E[\exp(-\eta e_t) F\biggl(\frac{1}{A_s}, s\leqq t\biggr) \biggr]
f(\eta) \frac{d\eta}{\eta},
\end{equation*}
we obtain \eqref{lpe:41}.

For Step 2, we replace $F \exp(\eta e_t)$ by $F$ in \eqref{lpe:42}.
Then, since $\exp(T_\alpha(B)_s) = e_s/(1+\alpha A_s),$ we obtain
\begin{equation*}
E[F(B_s,s\leqq t)] = E\biggl[ \exp(-\eta/e_t)
\exp\biggl(\frac{\eta e_t}{1+2\eta A_t/e_t}\biggr)
F(T_{2\eta/e_t}(B)_s, s\leqq t) \biggr]
\end{equation*}
and, setting $\eta=\alpha/2$,
\begin{equation} \label{lpe:43}
E[F(B_s,s\leqq t)] = E\biggl[
\exp\biggl( \frac{\alpha}{2} \biggl( \frac{e_t}{1+\alpha Z_t}
- \frac{1}{e_t} \biggr) \biggr)
F(T_{\alpha/e_t}(B)_s, s\leqq t) \biggr].
\end{equation}
\indent We now wish to find a function $\varphi:\R\times[0,\infty)\to\R$
such that
\begin{equation*}
E[\varphi(B_t,A_t) F(B_s, s\leqq t)] =
E[F(T_{\alpha/e_t}(B)_s, s\leqq t)].
\end{equation*}
With this aim in mind,
we replace $F$ by $\varphi F$ in \eqref{lpe:43}.
Then we have
\begin{equation*} \begin{split}
 & E[  \varphi(B_t,A_t) F(B_s, s\leqq t)] \\
 & = E\biggl[ \exp\biggl( \frac{\alpha}{2} \biggl( \frac{e_t}{1+\alpha Z_t}
 - \frac{1}{e_t} \biggr) \biggr)
 \varphi\biggl( T_{\alpha/e_t}(B)_t, \frac{A_t}{1+\alpha Z_t}\biggr)
 F(T_{\alpha/e_t}(B)_s, s\leqq t) \biggr].
\end{split} \end{equation*}
Hence, for our purpose, it is sufficient that
\begin{equation} \label{lpe:want}
\exp\biggl( \frac{\alpha}{2} \biggl( \frac{e_t}{1+\alpha Z_t}
 - \frac{1}{e_t} \biggr) \biggr)
\varphi\biggl( B_t-\log(1+\alpha Z_t), \frac{A_t}{1+\alpha Z_t}\biggr) =1
\end{equation}
and, if we take
\begin{equation*}
\varphi(b,u) = \exp\biggl( - \frac{\alpha}{2} \biggl( e^b -
\frac{1}{e^b+\alpha u} \biggr) \biggr),
\end{equation*}
we get \eqref{lpe:want}.   Therefore, we obtain
\begin{equation} \label{lpe:zero}
E[F(T_{\alpha/e_t}(B)_s,s \leqq t)] = E\biggl[ \exp\biggl(
- \frac{\alpha}{2} \biggl(e_t - \frac{1}{e_t+\alpha A_t} \biggr) \biggr)
F(B_s, s\leqq t) \biggr] ,
\end{equation}
which is precisely \eqref{lpe:rk} with $\mu=0$.

For Step 3, we again use the Cameron-Martin theorem to modify
the left hand side of \eqref{lpe:rk}:
\begin{align*}
 & E[ F(T_{\alpha/e_t^{(\mu)}}(B^{(\mu)})_s,s\leqq t)] \\
 & = E[F(T_{\alpha/e_t}(B)_s, s\leqq t) \exp(\mu B_t-\mu^2t/2)] \\
 & = e^{-\mu^2t/2} E[F(T_{\alpha/e_t}(B)_s, s\leqq t)
 \exp(\mu T_{\alpha/e_t}(B)_t) (1+\alpha Z_t)^\mu].
\end{align*}
We note $Z(T_{\alpha/e_t}(B))=Z(B)$ (Proposition \ref{lpp:easy})
and use \eqref{lpe:zero}.   Then we obtain
\begin{align*}
 & E[ F(T_{\alpha/e_t^{(\mu)}}(B^{(\mu)})_s,s\leqq t)] \\
 & = e^{-\mu^2t/2} E\biggl[F(B_s, s\leqq t) \exp(\mu B_t) (1+\alpha Z_t)^\mu
 \exp\biggl(-\frac{\alpha}{2}\biggl(e_t-\frac{1}{(1+\alpha Z_t)e_t}
 \biggr) \biggr) \biggr] \\
 & = E\biggl[ F(B_s^{(\mu)}, s\leqq t) (1+\alpha Z^{(\mu)}_t)^\mu
 \exp\biggl(-\frac{\alpha}{2}\biggl( e^{(\mu)}_t -
 \frac{1}{(1+\alpha Z_t^{(\mu)})e_t^{(\mu)}}
 \biggr) \biggr) \biggr].
\end{align*}
The proof is completed.
\end{proof}

We are now in a position to give a proof of the key proposition.

\medskip

\noindent{\it Proof of Proposition \ref{dp:key}}.\quad
At first we consider the case $\mu=0$.
We set $Q_t^{\omega,z}(\cdot)=P(\cdot|\mathcal{Z}_t, Z_t=z)$,
the regular conditional distribution given $\mathcal{Z}_t$.
Then, taking $F$ in \eqref{lpe:rk} as
$\varphi(1/A_t)G(Z_s,s\leqq t)$
for a non-negative Borel function $\varphi$ on $(0,\infty)$ and
for a non-negative functional $G$ in view of Proposition \ref{lpp:easy},
we obtain
\begin{equation*}
E^{Q_t^{\omega,z}}\biggl[ \exp(-\eta e_t) \varphi\biggl(\frac{1}{ze_t}
\biggr) \biggr] = E^{Q_t^{\omega,z}}\biggl[ \exp(-\eta/e_t)
\varphi\biggl(\frac{1}{ze_t}+\frac{2\eta}{e_t}\biggr) \biggr].
\end{equation*}
\indent Now we assume for simplicity that
the distribution of $e_t$ under $Q_t^{\omega,z}$ has a density $g_z(x)$
with respect to the Lebesgue measure.
Then we have
\begin{align*}
E^{Q_t^{\omega,z}}\biggl[ \exp(-\eta e_t) \varphi\biggl(\frac{1}{ze_t}
\biggr) \biggr] & = \int_0^\infty e^{-\eta x}
\varphi\biggl(\frac{1}{zx}\biggr) g_z(x) dx \\
 & = \int_0^\infty e^{-\eta/zx} \varphi(x)
 g_z\biggl(\frac{1}{zx}\biggr) \frac{1}{zx^2}dx
\end{align*}
and
\begin{align*}
 & E^{Q_t^{\omega,z}}\biggl[ \exp(-\eta/e_t)
\varphi\biggl(\frac{1}{ze_t}+\frac{2\eta}{e_t}\biggr) \biggr] \\
 & = \int_0^\infty e^{-\eta/x} \varphi\biggl(\frac{2\eta+1/z}{x}\biggr)
g_z(x) dx \\
 & = \int_0^\infty \exp\biggl(-\frac{\eta x}{2\eta+1/z}\biggr)
 \varphi(x) g_z\biggl(\frac{2\eta+1/z}{x}\biggr)
 \frac{2\eta+1/z}{x^2}dx.
\end{align*}
Since the function $\varphi$ is arbitrary, we obtain
\begin{equation*}
z^{-1} g_z(v/z) \exp(-\eta v/z) = (2\eta+1/z)
\exp\biggl(-\frac{\eta}{(2\eta+1/z)v}\biggr) g_z((2\eta+1/z)v),
\end{equation*}
where we have set $v=1/x$.
From the last identity, we obtain
\begin{equation*}
g_z(x) = \text{\rm const. } x^{-1}\exp\biggl(-\frac{1}{2z} \biggl(
x + \frac{1}{x} \biggr) \biggr)
\end{equation*}
by simple algebra and,
by using the integral representation \eqref{lpe:int-k}
for the Macdonald function, we obtain \eqref{de:gig} when $\mu=0$.

For a general value of $\mu$,
a standard argument with the Cameron-Martin theorem
leads us to the result. \hfill $\square$

\medskip

We finally show that
the semigroups of the diffusion processes $e^{(\mu)}$, $\xi^{(\mu)}$
and $Z^{(\mu)}$ satisfy some intertwining properties.
For some general discussions and examples of intertwinings
between Markov semigroups, see \cite{biane}, \cite{cpy}, \cite{dynkin},
\cite{rp} and \cite{yor-cr89} among others.

To present these results,
we denote by $I_z^{(\mu)}$ a generalized inverse Gaussian random variable
whose density is given by \eqref{de:gig} and
define the Markov kernels $\mathbb{K}_1^{(\mu)}$ and $\mathbb{K}_2^{(\mu)}$ by
\begin{equation*}
\mathbb{K}_1^{(\mu)}\varphi(z) = E[\varphi(I_z^{(\mu)})]
\quad \text{ and } \quad
\mathbb{K}_2^{(\mu)}\varphi(z) =
E\biggl[ \varphi\biggl( \frac{z}{I_z^{(\mu)}} \biggr) \biggr]
\end{equation*}
for a generic function $\varphi$ on $\R_{+}$.

\begin{thm}
Let $\{P_t^{(\mu)}\},$ $\{Q_t^{(\mu)}\}$ and $\{R_t^{(\mu)}\}$ be
the semigroup of the diffusion processes $e^{(\mu)},$ $\xi^{(\mu)}$
and $Z^{(\mu)},$ respectively{\rm.}
Then we have
\begin{equation*}
R_t^{(\mu)} \mathbb{K}_1^{(\mu)} = \mathbb{K}_1^{(\mu)} P_t^{(\mu)}
\quad \text{ and } \quad
R_t^{(\mu)} \mathbb{K}_2^{(\mu)} = \mathbb{K}_2^{(\mu)} Q_t^{(\mu)} .
\end{equation*}
\end{thm}

\begin{proof}
The key proposition (Proposition \ref{dp:key}) again plays an important role.
We give a proof for the first identity.
The second one is proven in a similar way.

For $s,t>0$,
we compute $E[\varphi(e_{s+t}^{(\mu)})|\mathcal{Z}_s^{(\mu)},
Z_s^{(\mu)}=z]$ in two ways.
First, we use the Markov property of $e^{(\mu)}$ to obtain
\begin{align*}
E[\varphi(e_{s+t}^{(\mu)})|\mathcal{Z}_s^{(\mu)}, Z_s^{(\mu)}=z] & =
E[E[\varphi(e_{s+t}^{(\mu)})|\Bo_s^{(\mu)}] |\mathcal{Z}_s^{(\mu)},
 Z_s^{(\mu)}=z] \\
 & = E[(P_t^{(\mu)}\varphi)(e_s^{(\mu)}) | \mathcal{Z}_s^{(\mu)},
 Z_s^{(\mu)}=z] \\
 & = E[(P_t^{(\mu)}\varphi)(I_z^{(\mu)})] \\
 & = (\mathbb{K}_1^{(\mu)} P_t^{(\mu)} \varphi)(z),
\end{align*}
where $\Bo_s^{(\mu)}=\sigma\{e_u^{(\mu)},u \leqq s\} =
\sigma\{B_u^{(\mu)},u\leqq s\}$ and
we have used Proposition \ref{dp:key} for the third equality.

Next, we note $E[\varphi(e_{s+t}^{(\mu)})|\mathcal{Z}_{s+t}^{(\mu)}] =
(\mathbb{K}_1^{(\mu)}\varphi)(Z_{s+t}^{(\mu)})$.
Then, by the Markov property of $Z^{(\mu)}$, we obtain
\begin{align*}
E[\varphi(e_{s+t}^{(\mu)})|\mathcal{Z}_s^{(\mu)}, Z_s^{(\mu)}=z] & =
E[\varphi(e_{s+t}^{(\mu)})|\mathcal{Z}_{s+t}^{(\mu)}] |
\mathcal{Z}_s^{(\mu)}, Z_s^{(\mu)}=z] \\
 & = E[(\mathbb{K}_1^{(\mu)}\varphi)(Z_{s+t}^{(\mu)}) |
 \mathcal{Z}_s^{(\mu)}, Z_s^{(\mu)}=z ] \\
 & = (R_t^{(\mu)} \mathbb{K}_1^{(\mu)} \varphi)(z)
\end{align*}
and the desired identity
$R_t^{(\mu)}\mathbb{K}_1^{(\mu)}=\mathbb{K}_1^{(\mu)}P_t^{(\mu)}$.
\end{proof}

\end{document}